\documentclass[a4paper,12pt] {article}
\usepackage{graphics}
\usepackage{graphicx}
\usepackage{pstricks}
\usepackage[english]{babel}
\usepackage{amsmath}
\usepackage{amssymb}
\usepackage{float}
\usepackage[caption=false]{subfig}
\usepackage{multirow}
\usepackage{bm}
\usepackage{url}
\usepackage{placeins}
\usepackage[abs]{overpic}
\usepackage{pict2e}
\usepackage{pgfplots}

\definecolor{lgray}{gray}{0.85}
\definecolor{dgray}{RGB}{128,128,128}
\definecolor{lblue}{RGB}{0,191,191}
\definecolor{lblue2}{RGB}{77,191,237}
\definecolor{dblue}{RGB}{0,114,189}
\definecolor{orange}{RGB}{222,125,0}
\definecolor{dorange}{RGB}{217,83,25}
\definecolor{dred}{RGB}{189,18,35}
\definecolor{dgreen}{RGB}{119,172,48}
\definecolor{dyellow}{RGB}{255,215,0}

\pgfplotsset{scaled ticks = false}

\pgfplotsset{%
    ,compat=1.11
    ,colormap={parula}{%
        rgb=(0.2081,0.1663,0.5292)rgb=(0.2116,0.1898,0.5777)rgb=(0.2123,0.2138,0.627)
        rgb=(0.2081,0.2386,0.6771)rgb=(0.1959,0.2645,0.7279)rgb=(0.1707,0.2919,0.7792)
        rgb=(0.1253,0.3242,0.8303)rgb=(0.0591,0.3598,0.8683)rgb=(0.0117,0.3875,0.882)
        rgb=(0.006,0.4086,0.8828) rgb=(0.0165,0.4266,0.8786)rgb=(0.0329,0.443,0.872)
        rgb=(0.0498,0.4586,0.8641)rgb=(0.0629,0.4737,0.8554)rgb=(0.0723,0.4887,0.8467)
        rgb=(0.0779,0.504,0.8384) rgb=(0.0793,0.52,0.8312)  rgb=(0.0749,0.5375,0.8263)
        rgb=(0.0641,0.557,0.824)  rgb=(0.0488,0.5772,0.8228)rgb=(0.0343,0.5966,0.8199)
        rgb=(0.0265,0.6137,0.8135)rgb=(0.0239,0.6287,0.8038)rgb=(0.0231,0.6418,0.7913)
        rgb=(0.0228,0.6535,0.7768)rgb=(0.0267,0.6642,0.7607)rgb=(0.0384,0.6743,0.7436)
        rgb=(0.059,0.6838,0.7254) rgb=(0.0843,0.6928,0.7062)rgb=(0.1133,0.7015,0.6859)
        rgb=(0.1453,0.7098,0.6646)rgb=(0.1801,0.7177,0.6424)rgb=(0.2178,0.725,0.6193)
        rgb=(0.2586,0.7317,0.5954)rgb=(0.3022,0.7376,0.5712)rgb=(0.3482,0.7424,0.5473)
        rgb=(0.3953,0.7459,0.5244)rgb=(0.442,0.7481,0.5033) rgb=(0.4871,0.7491,0.484)
        rgb=(0.53,0.7491,0.4661)  rgb=(0.5709,0.7485,0.4494)rgb=(0.6099,0.7473,0.4337)
        rgb=(0.6473,0.7456,0.4188)rgb=(0.6834,0.7435,0.4044)rgb=(0.7184,0.7411,0.3905)
        rgb=(0.7525,0.7384,0.3768)rgb=(0.7858,0.7356,0.3633)rgb=(0.8185,0.7327,0.3498)
        rgb=(0.8507,0.7299,0.336) rgb=(0.8824,0.7274,0.3217)rgb=(0.9139,0.7258,0.3063)
        rgb=(0.945,0.7261,0.2886) rgb=(0.9739,0.7314,0.2666)rgb=(0.9938,0.7455,0.2403)
        rgb=(0.999,0.7653,0.2164) rgb=(0.9955,0.7861,0.1967)rgb=(0.988,0.8066,0.1794)
        rgb=(0.9789,0.8271,0.1633)rgb=(0.9697,0.8481,0.1475)rgb=(0.9626,0.8705,0.1309)
        rgb=(0.9589,0.8949,0.1132)rgb=(0.9598,0.9218,0.0948)rgb=(0.9661,0.9514,0.0755)
        rgb=(0.9763,0.9831,0.0538)
        }
    }
    
\pgfplotsset{compat=1.3}
\begin{document}
\title{\Large\textbf{Numerical Continuation in Nonlinear Experiments using Local Gaussian Process Regression}}
\author{\textbf{L.~Renson$^{1,}$\footnote{Corresponding author: l.renson@bristol.ac.uk}, J.~Sieber$^{2}$, D.A.W.~Barton$^{1}$, A.D.~Shaw$^{3}$, S.A.~Neild$^{1}$} \vspace{2mm}\\
\normalsize{$^{1}$Faculty of Engineering, University of Bristol, UK.}\\
\normalsize{$^{2}$Centre for Systems, Dynamics and Control, College of Engineering, University of Exeter, UK.}\\
\normalsize{$^{3}$College of Engineering, Swansea University, UK.}}
\date{}
\maketitle

\begin{abstract}
\noindent
Control-based continuation (CBC) is a general and systematic method to probe the dynamics of nonlinear experiments. In this paper, CBC is combined with a novel continuation algorithm that is robust to experimental noise and enables the tracking of geometric features of the response surface such as folds. The method uses Gaussian process regression to create a local model of the response surface on which standard numerical continuation algorithms can be applied. The local model evolves as continuation explores the experimental parameter space, exploiting previously captured data to actively select the next data points to collect such that they maximise the potential information gain about the feature of interest. The method is demonstrated experimentally on a nonlinear structure featuring harmonically-coupled modes. Fold points present in the response surface of the system are followed and reveal the presence of an isola, i.e. a branch of periodic responses detached from the main resonance peak.\\

\end{abstract}

\section{Introduction}
Numerical continuation is a popular and well-esta\-blished method to systematically investigate the behaviour of nonlinear dynamical systems and perform bifurcation analysis~\cite{KuznetsovBook,SeydelBook}. At a basic level, numerical continuation finds the solutions of a zero-problem $\mathbf{f}(\mathbf{x},\lambda) = \mathbf{0}$ where $\mathbf{x}$ are the system states, and tracks the evolution of the solutions as the parameter $\lambda$ is varied. Based on a mathematical model, the long-term behaviours of a system, such as steady-states and periodic responses, can easily be represented by such a zero-problem. Bifurcations, and hence stability changes, in those behaviours can then be detected along the solution path and in turn tracked by adding suitable constraint condition(s) and free parameter(s). The principles of numerical continuation are extremely general such that the method has been applied to a wide range of problems across engineering and the applied sciences as, for instance, in bio-chemistry~\cite{Godwin17}, physics~\cite{Krauskopf05}, mechanics~\cite{Peeters09} and fluid dynamics~\cite{Huntley17}.\\ 

Without the need for a mathematical model, control-based continuation (CBC) is a means to define a zero-problem based on the inputs and outputs of an experiment, thereby allowing the principles of numerical continuation to be applied to a physical system directly during experimental tests. The fundamental idea of CBC is to use feedback control to stabilise the dynamics of the experiment whilst making the control system \textit{non-invasive} such that it does not modify the position in parameter space of the responses of the open-loop experiment of interest. This non-invasiveness requirement defines a zero-problem whose solutions can be found and tracked in the experiment using the same path-following principles and methods as in the numerical context.\\

CBC is similar in principle to other methods such as the famous OGY (Ott, Grebogi, Yorke) control technique that has been extensively used to stabilise unstable periodic responses embedded in chaotic attractors~\cite{Ott90}. The OGY method was coupled to continuation algorithms in~\cite{Misra08}, but the application of this method to general nonlinear experiments\footnote{The term ``\textit{nonlinear experiment}'' refers to an experiment for which any model describing its behaviour has to be nonlinear to be consistent with observations.} remains challenging due to the particular form of control used in the OGY technique. Other examples of control techniques used to experimentally measure unstable responses are the Pyragas delayed feedback control~\cite{Pyragas06} and phase-locked loops~\cite{PLLBook}. The latter has recently been applied to several nonlinear mechanical systems~\cite{MojrzischPAMM,Peter17,Denis18}. Although these methods and CBC share a number of similarities, CBC does not assume any particular form of control.\\ 

CBC was first proposed by Sieber and Krauskopf~\cite{Sieber08}, and experimentally demonstrated on a parametrically-excited pendulum~\cite{Sieber10}. The method has since been successfully applied to a range of mechanical systems, including an impact oscillator~\cite{Bureau13,Bureau14,Elmegard14}, oscillators with magnetic nonlinearities~\cite{Barton13,Renson17} and a cantilever beam with a nonlinear mechanism at its free tip~\cite{Renson19}. Through those studies, CBC proved to be a versatile technique that can extract important nonlinear dynamic features such as nonlinear frequency response curves~\cite{Barton10}, nonlinear normal modes~\cite{Renson16,RensonIMAC2016} and isola~\cite{Renson19} directly in the experiment. However, the systematic application of CBC to general nonlinear experiments remains challenging. Most existing continuation algorithms are ideal only in a numerical context where the solution path is smooth and derivatives can be evaluated to high precision. This is not easily achievable in experiments where solutions and derivative estimates are corrupted by measurement noise. Schilder \textit{et al.}~\cite{Schilder15} discussed the effect of noise on continuation algorithms. In particular, the tangential prediction and orthogonal correction steps of the commonly-used pseudo-arclength continuation algorithm were shown to perform poorly in a noisy experimental context. Similarly, continuation step size control techniques that reduce continuation steps when convergence is not achieved were also shown to be inadequate as reducing the step size usually makes noise distortions even more apparent. Schilder \textit{et al.} proposed alternative numerical strategies that are more robust to noise~\cite{Schilder15}. These new strategies significantly improve the robustness of CBC as they are able to trace solution paths that are no longer smooth due to noise. These strategies work well for low levels of noise but are not sufficiently robust to find and track dynamic features that are very sensitive to noise as, for instance, bifurcations which are defined in terms of derivatives.\\

This paper proposes an algorithm that not only makes CBC more robust to noise but also enables the tracking of general dynamic features that are not directly measurable in the experiment due to perturbations from noise for instance. The proposed approach is fundamentally different from the one taken by Schilder and co-workers. Multivariate regression techniques are exploited to locally model the response surface of the experiment, providing local models that are smooth and cheap to evaluate and thus enabling the use of established numerical continuation techniques. The proposed method is inspired by the work reported in~\cite{Renson17} where a single-parameter cubic polynomial regression was exploited to capture the geometry of the response surface at a fixed forcing frequency and used to detect and track a limit-point bifurcation in parameter space. The algorithm presented in~\cite{Renson17} is however limited to bifurcation curves with simple geometries. For instance, the algorithm cannot follow the bifurcation through a cusp. The method proposed here overcomes this limitation.\\

A key difficulty in using online models is to devise a strategy to collect the experimental data necessary to build them. We use Gaussian Process regression (GPR) techniques to address this challenge. GPR has many desirable features, such as the ease of extension to models with multiple inputs and outputs, the ease of expressing uncertainty, the ability to capture a wide range of behaviours using a simple (hyper-)parameterisation, and a natural Bayesian interpretation. Here, based on the data points already captured, GPR will allow us to determine where to collect new data points to maximise the potential information they will provide about the dynamic feature of interest. This active selection of the data based on our current knowledge and the feature of interest contrasts with the approaches currently found in the literature~\cite{MSSPrev,NoelReview} where data collection and identification are two activities often performed separately.\\

The proposed algorithm is presented in Section~\ref{sec:theory} and demonstrated on a nonlinear mechanical structure composed of a cantilever beam with a nonlinearity attached at its free end (Sections~\ref{sec:rig} and~\ref{sec:results}). The natural frequencies of the first two modes of the structure are almost in a 3:1 ratio, which leads to strong harmonic couplings between these modes and the presence of complicated nonlinear dynamic behaviours. In particular, a branch of stable periodic responses detached from the main resonance peak can be observed. Our new algorithm is employed to track fold points present in the response surface of the system. The obtained curves are found to have a complex geometry due to the presence of the modal interaction and to reveal the presence of the isola.\\

\section{Tracking dynamic features using an online regression-based continuation algorithm}\label{sec:theory}
The zero-problem defining the responses tracked in the experiment is now explained (Section~\ref{sec:theory_cbc}). The algorithm used to experimentally solve and continue the solution of this zero-problem is then presented in Section~\ref{sec:algo_overview}. Detailed discussions of the important components of this algorithm are given in Sections~\ref{sec:algo_gpr} -- \ref{sec:algo_improve}.

\subsection{Definition of the zero-problem}\label{sec:theory_cbc}
\begin{figure}[tbp]
\centering
\includegraphics[width=0.45\textwidth]{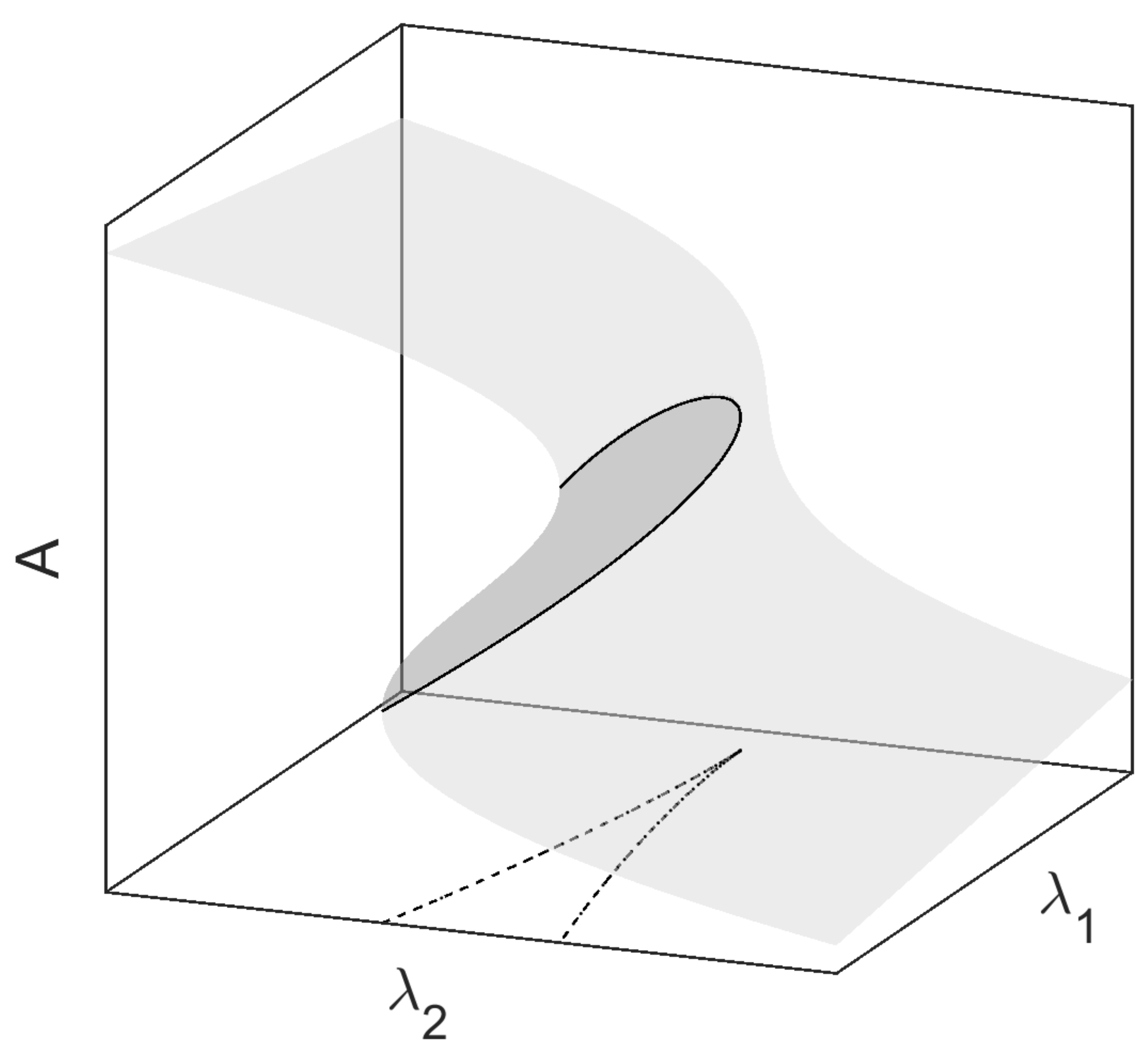}
\caption{Typical response surface that can be found in the experiment. $A$ is a measure of the response of the system and $\lambda_1$ and $\lambda_2$ are parameters, here selected as $\lambda_1 = \omega$ and $\lambda_2 = \Gamma$.  Fold curve (\textcolor{black}{$\boldsymbol{-}$}) and its projection (\textcolor{black}{$\boldsymbol{--}$}) in parameter space.}
\label{fig:Surface_th}
\end{figure}

The response surface of an experiment such as the one considered here  presents the generic form shown in Figure~\ref{fig:Surface_th}, where $A$ is the system response amplitude and $(\lambda_1, \; \lambda_2)$ are parameters. Such response surface can be easily obtained experimentally using an established CBC algorithm~\cite{Renson17}. In this paper, we are interested in directly tracking geometric features of that response surface during the experiment and, in particular, the fold curve represented in solid black (\textcolor{black}{$\boldsymbol{-}$}). Responses that lie on that curve satisfy the scalar constraint equation
\begin{equation}
\frac{d \lambda_2 (\lambda_1, A)}{dA} = 0.
\label{eq:LP_cond}
\end{equation}
The response amplitude $A$ is a natural choice for the parameterisation as it defines the fold curve uniquely in many experiments (as in Figure~\ref{fig:Surface_th}). The parameter $\lambda_1$ serves as the free parameter for the continuation and, in the present experimental context, corresponds to the forcing frequency $\omega$. $\lambda_2$ is the external harmonic force excitation amplitude $\Gamma$.\\

\setlength{\unitlength}{1cm}
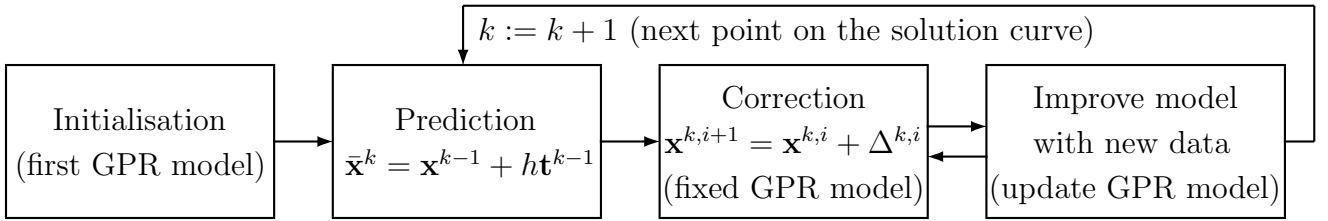
\begin{figure*}[tpb]
    \centering
    \begin{picture}(17.2, 2.9)
    \thicklines
    \put(0, 0){\framebox(3.5,2)}
    \put(0.6, 1.15){Initialisation}
    \put(0.1, 0.55){(first GPR model)}
    \put(3.5, 1){\vector(1, 0){0.8}}
    \put(4.3, 0){\framebox(3.5,2)}
    \put(5.1, 1.15){Prediction}
    \put(4.45, 0.55){$\bar{\mathbf{x}}^k = \mathbf{x}^{k-1} + h \mathbf{t}^{k-1}$}
    \put(7.8, 1){\vector(1, 0){0.8}}
    \put(8.6, 0){\framebox(3.5,2)}
    \put(9.4, 1.45){Correction}
    \put(8.65, 0.85){$\mathbf{x}^{k,i+1} = \mathbf{x}^{k,i} + \Delta^{k,i}$}
    \put(8.62, 0.25){(fixed GPR model)}
    \put(12.1, 1.2){\vector(1, 0){0.8}}
    \put(12.9, 0.8){\vector(-1, 0){0.8}}
    \put(12.9, 0){\framebox(3.9,2)}
    \put(13.5, 1.45){Improve model}
    \put(13.5, 0.85){with new data}
    \put(12.92, 0.25){(update GPR model)}
    \put(16.8, 1){\line(1,0){0.4}}
    \put(17.2, 1){\line(0,1){1.8}}
    \put(17.2, 2.8){\line(-1,0){11.2}}
    \put(6.0, 2.8){\vector(0,-1){0.8}}
    \put(6.2, 2.35){$k:=k+1$ (next point on the solution curve)}
    \end{picture}
    \caption{Overall structure of the regression-based continuation algorithm proposed here. $\mathbf{x}^k$ is the k$^{th}$ solution found using the online regression-based continuation, $h$ is the continuation step size, $\mathbf{t}^k$ is the tangent to the solution curve at $\mathbf{x}^k$. $\Delta^{k,i}$ is the correction applied to $\mathbf{x}^{k,i}$ at the i$^{th}$ iteration of the correction algorithm.}
    \label{fig:algo}
\end{figure*}

The choice of $\omega$ and $A$ as independent variables in Eq.~\eqref{eq:LP_cond} also stems from the nature of these variables in the experiment. More precisely, the response amplitude $A$ is indirectly imposed by the amplitude of the reference signal of the control system used in CBC (in particular, one of its fundamental harmonic components). As such, $A$ and $\omega$ can both be viewed as `controllable inputs' to the experiment. In contrast to continuation applied to a mathematical model, the force amplitude $\Gamma$ is here considered as a quantity that is difficult to set directly and hence one that is measured from the experiment rather than imposed.

\subsection{Overview of the online regression-based algorithm}\label{sec:algo_overview}
Finding data points experimentally that satisfy Eq.~\eqref{eq:LP_cond} can be difficult and error prone due to the presence of noise affecting derivative calculations. This issue is addressed here by creating, online (i.e. while the experiment is running), regression models that capture the local dependence of the applied force amplitude as a function of the forcing frequency and response amplitude. Derivatives can be effectively and accurately calculated for these models, which in turn allows us to find and then track the fold curve in the experiment using standard numerical continuation algorithms.\\

The principal steps of the proposed algorithm are shown in Figure~\ref{fig:algo}. The algorithm is initialised by collecting a user-defined number of data points $n_0$ distributed in a regular pattern around a starting point which is assumed to be close to a fold. This first set of experimental data is then used to create the first GPR model and estimate the regression hyper-parameters (Section~\ref{sec:algo_gpr}). Using this model, a fold point is also found by solving Eq.~\eqref{eq:LP_cond} using a standard Newton-like algorithm. From this first point, standard predictor-corrector continuation algorithms can be exploited (Section~\ref{sec:algo_cont}). However, unlike in numerical simulations, the correction step is followed by a data collection step that aims to make the solution of the zero-problem~\eqref{eq:LP_cond} robust to new data points (Section~\ref{sec:algo_improve}). The addition of each new data point to the local GPR model is followed by a correction procedure to update the solution of Eq.~\eqref{eq:LP_cond}. When no additional data is needed, a new prediction step is performed to progress along the solution curve.\\

\subsection{Gaussian Process Regression}\label{sec:algo_gpr}
The amplitude of excitation $\Gamma$ is locally modelled as a function of the response amplitude and excitation frequency using GPR. A Gaussian process is a probabilistic model that can be used to capture a wide range of nonlinear functions from input-output data without any explicit assumptions on their mathematical relationship~\cite{RasmussenBook}. The amplitude of excitation, $\Gamma$, is modeled by the distribution\footnote{Rigorously, our notations should distinguish the model output from the exact force amplitude $\Gamma$ as the model is only an approximation of the truth. However, in the present experimental context, the exact force amplitude is unavailable to us. As such, to keep our notations simple, the model output will also be denoted $\Gamma$.} 
\begin{equation}
    \Gamma(\mathbf{x}) \sim \mathcal{GP}\left(m(\mathbf{x}), \kappa (\mathbf{x},\mathbf{x}')\right),
    \label{eq:GP_Gamma}
\end{equation}
where $\mathbf{x} = \left( \omega, A \right)$ is the vector of inputs, $m(\mathbf{x})$ and $\kappa(\mathbf{x},\mathbf{x})$ are the mean and covariance functions. For the data point $i$, the measured excitation amplitude, denoted $F_i$, is assumed to differ from the function values $\Gamma_i$ by the presence of an additive Gaussian noise with zero mean and variance $\sigma_n^2$ such that
\begin{equation}
F_i = \Gamma_i(\mathbf{x}) + \epsilon_i \quad \text{with} \quad \epsilon_i \sim \mathcal{N}\left(0, \sigma_n^2 \right).
    \label{eq:noise_model}
\end{equation}
This `noise' captures in fact measurement errors, which not only include measurement noise but also other inaccuracies that arise in the collection of data points using CBC (see Section~\ref{sec:results}).\\

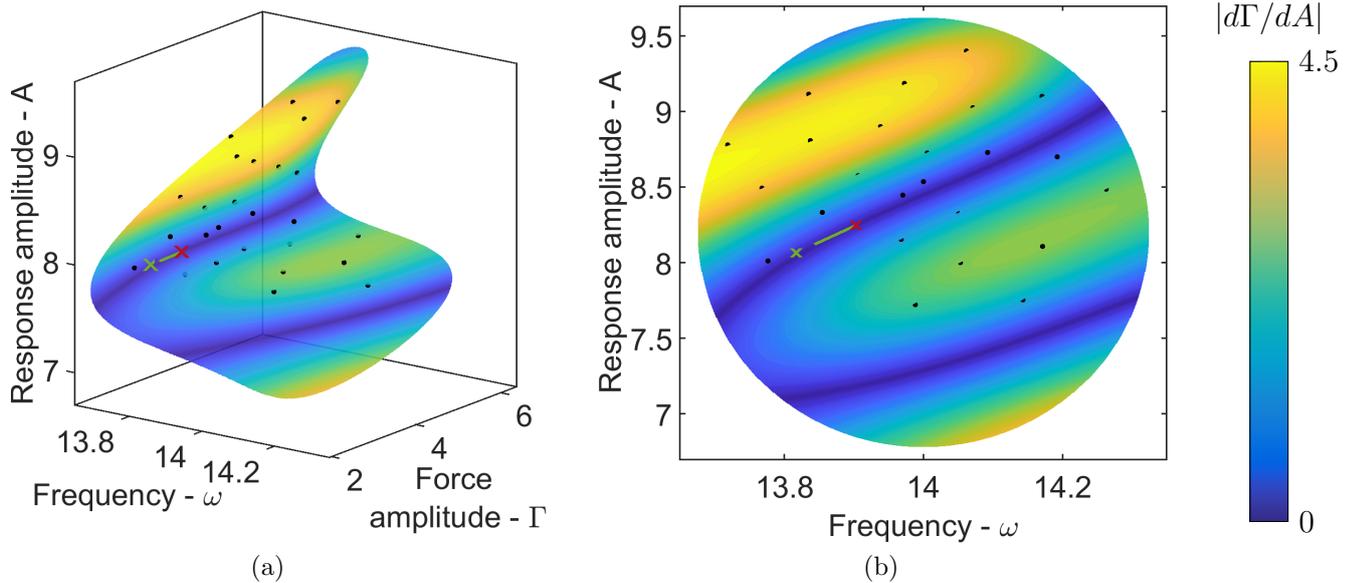
\begin{figure*}[t]
\centering
\begin{tabular*}{1.\textwidth}{@{\extracolsep{\fill}} c c c}
\subfloat[]{\label{ellipse}\includegraphics[width=0.45\textwidth]{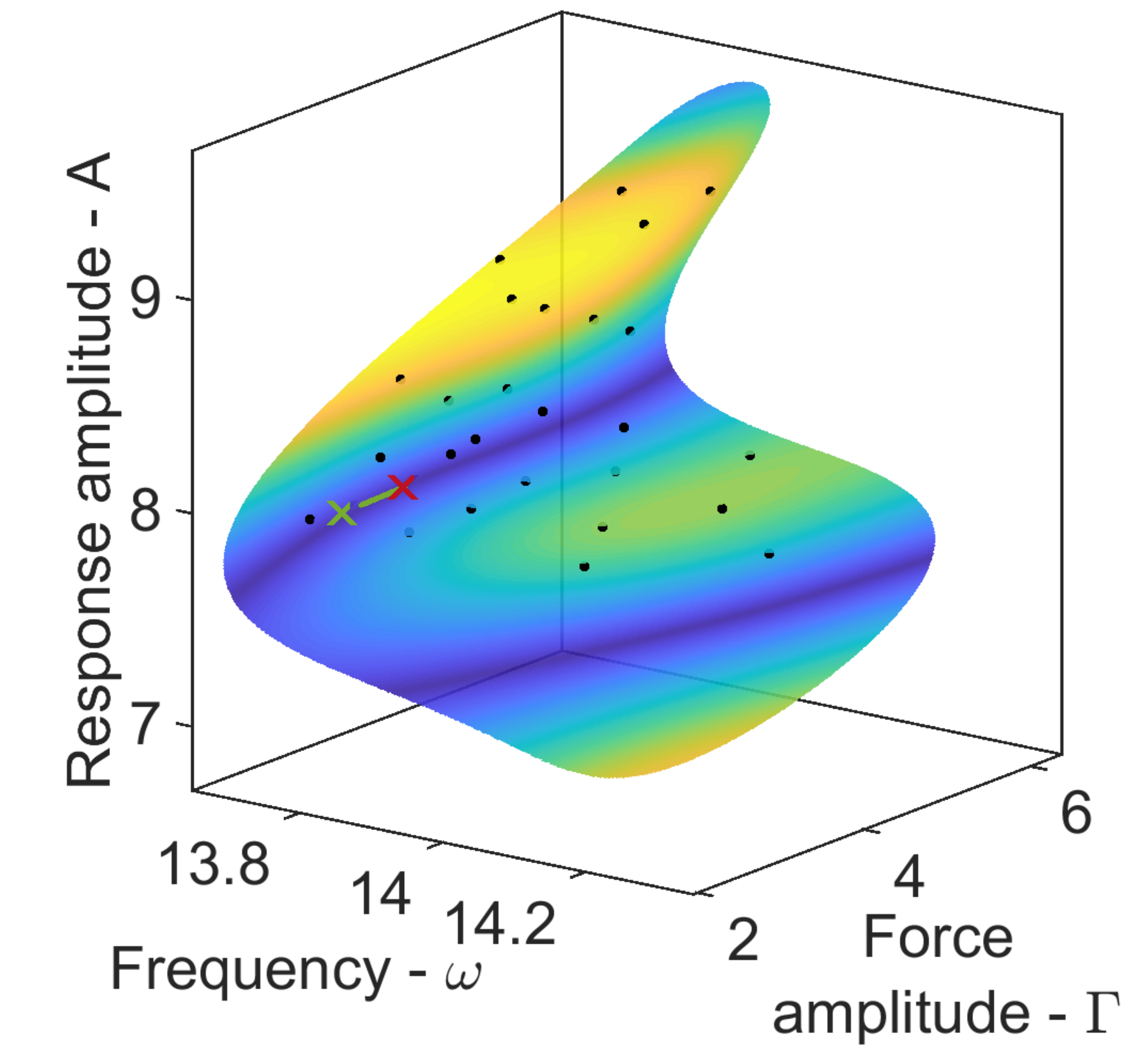}} &
\subfloat[]{\label{var}\includegraphics[width=0.45\textwidth]{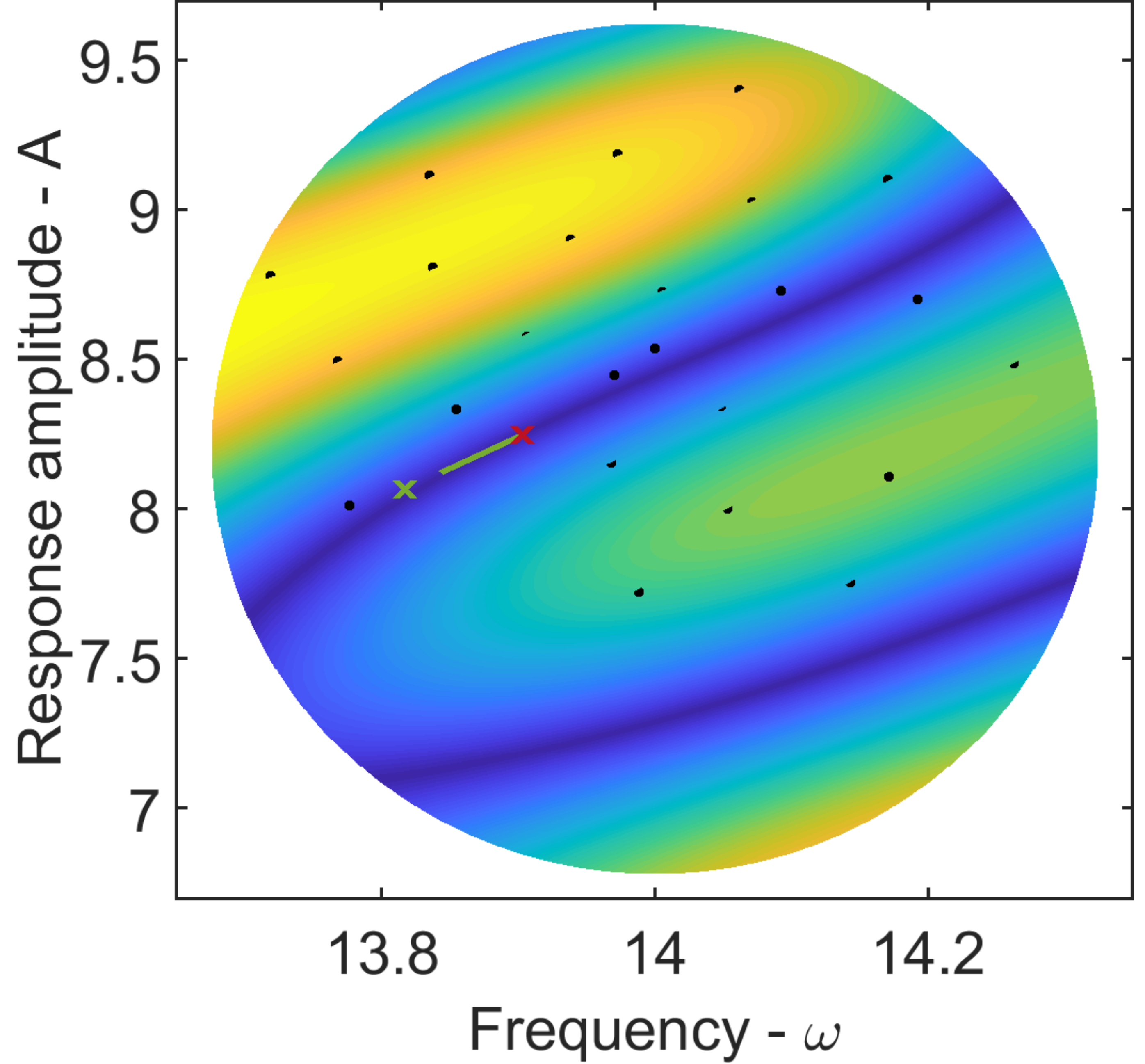}} & 
\begin{tikzpicture}
\begin{axis}[
    hide axis,
    scale only axis,
    height=0pt,
    width=0pt,
    colorbar,
    colorbar style={height = 6.1cm, ytick={0, 4.5}, title={$|d\Gamma / dA|$}},
    point meta min=0,
    point meta max=4.5
    ]
    \addplot [draw=none] coordinates {(0,0)};
\end{axis}
\end{tikzpicture}
\end{tabular*}
\caption{Illustration of the local model obtained using GPR at the initialisation of the algorithm. (a) and (b) correspond to different views of the same model. The colour map represents the absolute value of Eq.~\eqref{eq:LP_cond}, i.e. $|d\Gamma / dA|$. (\textcolor{black}{$\boldsymbol{\bullet}$}) Data points used to create the model. (\textcolor{red}{$\boldsymbol{\times}$}) Fold point. (\textcolor{dgreen}{$\boldsymbol{\times}$}) Predicted fold point. (\textcolor{dgreen}{$\boldsymbol{-}$}) Tangent vector to the solution curve at fold point (\textcolor{red}{$\boldsymbol{\times}$}).}
\label{fig:continuation}
\end{figure*}

Given a data set $\mathcal{D}$ with $n$ measurements
\begin{equation*}
    \left\{ (\mathbf{x}_1,F_1), (\mathbf{x}_2,F_2), ..., (\mathbf{x}_n,F_n) \right\},
\end{equation*}
grouped in an input matrix $\mathbf{X}=[\mathbf{x}_1, \mathbf{x}_2, ..., \mathbf{x}_n]$ and an output vector $\mathbf{F} = [F_1, F_2, ..., F_n]^T $, the prediction of the force amplitude at $n_{\star}$ unmeasured inputs $\mathbf{X}_{\star}$ is given by the mean of the predictive distribution as
\begin{equation}
\boldsymbol{\Gamma}_{\star} = \boldsymbol{\kappa} \left( \mathbf{X}, \mathbf{X}_{\star}\right)^T  \boldsymbol{\kappa} \left( \mathbf{X}, \mathbf{X} \right) ^{-1} \mathbf{F}. 
    \label{eq:mean_gamma} 
\end{equation}
where $(\cdot)^T$ represents the transpose operation. $\boldsymbol{\kappa} \left( \mathbf{X}, \mathbf{X}_{\star}\right)$ corresponds to the $n\times n_{\star}$ matrix resulting from the application of the covariance function $\kappa \left( \mathbf{x}_i, \mathbf{x}_{\star j}\right)$ for all $i$ and $j$. Similarly for $\boldsymbol{\kappa} \left( \mathbf{X}, \mathbf{X}\right)$. Note that GPR differs from parameter estimation techniques and the regression method used in~\cite{Renson17} because measured data points are needed to perform predictions (see Eq.~\eqref{eq:mean_gamma}).\\

An effective computation of Eq.~\eqref{eq:mean_gamma} can be achieved using the Cholesky factorization to decompose the covariance matrix $\boldsymbol{\kappa} \left( \mathbf{X}, \mathbf{X} \right)$ into a lower-triangular matrix and its conjugate transpose~\cite{RasmussenBook}. The former is then stored and used to efficiently compute the inverse for different unmeasured inputs $\mathbf{x}_{\star}$. The Cholesky decomposition can also be efficiently updated when data points are added or removed from the data set $\mathcal{D}$ as in Section~\ref{sec:algo_improve}. Note that the inputs $(\omega, A)$ of the GPR model are assumed to be noise-free but they are in reality measured quantities corrupted by noise. It is possible to extend the GPR model used here to address noisy inputs~\cite{McHutchon2011}; however this is not considered here as it was unnecessary for the system considered.\\

Although GPR is a non-parametric approach that does not make assumptions on the functional form of the modelled function, assumptions regarding its smoothness are introduced in the covariance function and the choice of hyper-parameters $\boldsymbol{\theta}$. The covariance function considered here is the widely-used, infinitely-differentiable squared-exponential (SE) or radial basis function (see, for instance, Eq.~(5.1)~in~\cite{RasmussenBook}). With two-dimensional inputs and observation noise, this covariance function includes four hyper-parameters $\boldsymbol{\theta} = (\sigma_n^2, \sigma_f^2, l_{\omega}, l_A)$. These parameters have a clear physical interpretation. In particular, $l_{\omega}$ and $l_A$ represent characteristic length scales that express the distance to be travelled along a particular axis in the input space for covariance values to become uncorrelated. $\sigma_n^2$ is the variance associated with measurement errors, and $\sigma_f^2$ is a magnitude factor. The hyper-parameters $\boldsymbol{\theta}$ of the models were obtained by maximizing the $\text{log}$ marginal likelihood
\begin{equation}
\log p( \mathbf{F} | \mathbf{X}, \theta) = \frac{1}{2} \mathbf{F}^T \boldsymbol{\kappa} \left( \mathbf{X}, \mathbf{X} \right)^{-1} \mathbf{F}   -\frac{1}{2} \log |\boldsymbol{\kappa} \left( \mathbf{X}, \mathbf{X} \right)| - \frac{n}{2} \log 2 \pi
\label{eq:log_marg}
\end{equation}
\noindent where the first term represents the data fit of the model and the second term penalizes the complexity of the model~\cite{RasmussenBook}. The optimization was performed using a Quasi-Newton algorithm. Hyper-parameters were determined at the start of the algorithm using the first set of $n_0$ data points captured. They were then kept constant during a continuation run. We note that Gaussian priors on the hyper-parameters (hyper-priors) can be used. However, with limited knowledge of the actual values and hence large covariances they had limited influence on the optimal $\boldsymbol{\theta}$ values.\\

As an example of a GPR model, Figure~\ref{fig:continuation}(a) shows a local model obtained at the start of a continuation run after collecting 25 data points (\textcolor{black}{$\boldsymbol{\bullet}$}) regularly distributed in input space $(\omega, A)$ around the estimated location of a fold point. A projection of the model in the two-dimensional input space is also given in Figure~\ref{fig:continuation}(b). The surface is coloured according to $|d\Gamma / dA|$ to highlight the regions where Eq.~\eqref{eq:LP_cond} is satisfied. According to the model, fold points are expected in two distinct regions (in dark-blue). The one where the response amplitude is higher is where the actual folds are located. The lower one is in fact located outside the data set and is an artifact feature created by the regression. These artifacts do not affect the algorithm as long as the continuation steps are small enough to stay within the available data sets. In Section~\ref{sec:algo_improve}, we will discuss how to improve such GPR models with additional experimental data but first we discuss the continuation approach.\\

\subsection{Numerical continuation}\label{sec:algo_cont}
The continuation problem is to solve and track the solutions of Eq.~\eqref{eq:LP_cond}. Starting from a known solution $\mathbf{x}^{k-1}$, the next point along the solution path is predicted to be
\begin{equation}
    \tilde{\mathbf{x}}^k = \mathbf{x}^{k-1} + h \; \mathbf{t}^{k-1}
    \label{eq:prediction}
\end{equation}
where $h$ is the continuation step size and $\mathbf{t}^{k-1}$ is the tangent vector to the solution curve at $\mathbf{x}^{k-1}$. The prediction will not in general satisfy the zero-problem such that the prediction must be corrected using a Newton-like algorithm. However, to apply Newton iterations, an extra equation has to be added to Eq.~\eqref{eq:LP_cond}. The equation used here is the so-called pseudo-arclength condition
\begin{equation}
    \mathbf{t}^{k-1} . (\mathbf{x} - \mathbf{x}^{k-1}) - h = 0,
    \label{eq:psa}
\end{equation}
which constraints the corrections made to $\tilde{\mathbf{x}}^k$ to be perpendicular to $\mathbf{t}^{k-1}$.\\

Figure~\ref{fig:continuation} illustrates a fold point (\textcolor{dgreen}{$\boldsymbol{\times}$}) predicted from a previously computed point (\textcolor{dred}{$\boldsymbol{\times}$}). The tangent vector used for this prediction is shown in green (\textcolor{dgreen}{$\boldsymbol{-}$}). It was found that GPR models could diverge quickly from the actual response surface outside the set of collected data points (see Section~\ref{sec:algo_gpr}). As such, the size of the continuation steps were taken such that the prediction and correction iterations remained within the cloud of available data points. Note that it is possible to allow larger continuation steps that leave the current data set by collecting new data points after the prediction step, although this is not used here. After finding a solution to Eqs.~\eqref{eq:LP_cond} and~\eqref{eq:psa}, additional data points are collected to refine this solution --- this is now discussed. 

\subsection{Improve solution with new data}\label{sec:algo_improve}
To make sure that the solution of Eq.~\eqref{eq:LP_cond} is not an artifact of the model and does not critically depend on the current data, new data points are collected after the correction step of the continuation algorithm. To this end, a fixed number $n_{test}$ of prospective data points, $\mathbf{x}_{\star}^i$ with $i=1, ..., n_{test}$, uniformly distributed in an ellipse around the current solution of the continuation problem are considered for data collection. The principal axis of this elliptical domain were chosen equal to twice the length scale hyper-parameters $l_{\omega}$ and $l_A$. The data point that most influences the results is deemed to be the most interesting point and is then experimentally collected using the established CBC technique summarised in Section~\ref{sec:results_basic}.\\

To determine the sensitivity of the zero-problem~\eqref{eq:LP_cond} to new experimental data, an artificial measurement of the force amplitude, $\bar{\mathbf{F}_i}$, is created for each candidate data point $\mathbf{x}_{\star}^i$ using
\begin{equation}
   \bar{F}^i = \Gamma^i_{\star}(\mathbf{x}_{\star}^i) + \sqrt{\text{var}[\Gamma^i_{\star}(\mathbf{x}_{\star}^i)]},
    \label{eq:artificial_meas}
\end{equation}
where $\Gamma^i_{\star}$ depends on currently available data and is given by Eq.~\eqref{eq:mean_gamma}. The variance at a particular input point $\mathbf{x}_{\star}^i$ is given by
\begin{equation}
    \text{var}[\Gamma_{\star}^i] = \kappa \left( \mathbf{x}_{\star}^i , \mathbf{x}_{\star}^i \right) - \boldsymbol{\kappa} \left( \mathbf{X}, \mathbf{x}_{\star}^i\right)^T  \boldsymbol{\kappa} \left( \mathbf{X}, \mathbf{X} \right) ^{-1} \boldsymbol{\kappa} \left( \mathbf{X}, \mathbf{x}_{\star}^i\right),
    \label{eq:var}
\end{equation}
which is independent of any previous and future measurement of the force amplitude. The artificial measurements are individually added to the GPR model and their effect on the zero-problem assessed using
\begin{equation}
\beta = \left|\frac{d \tilde{\Gamma}}{dA}\left(\mathbf{x}^k\right) - \frac{d \Gamma}{dA}\left(\mathbf{x}^k\right)\right|
\label{eq:comp}    
\end{equation}
where $\mathbf{x}^k$ is the current solution of the continuation problem, $\tilde{\Gamma}$ is the GPR model including the artificial measurement and $\Gamma$ is the GPR model without (i.e. containing only the experimental data). $\beta$ is the zero-problem sensitivity to new data and is directly used to assess the potential information gained by each artificial measurement.\\
\begin{figure*}[ht]
\centering
\begin{tabular*}{1.\textwidth}{@{\extracolsep{\fill}} c c c}
\subfloat[]{\label{sens1}\includegraphics[width=0.32\textwidth]{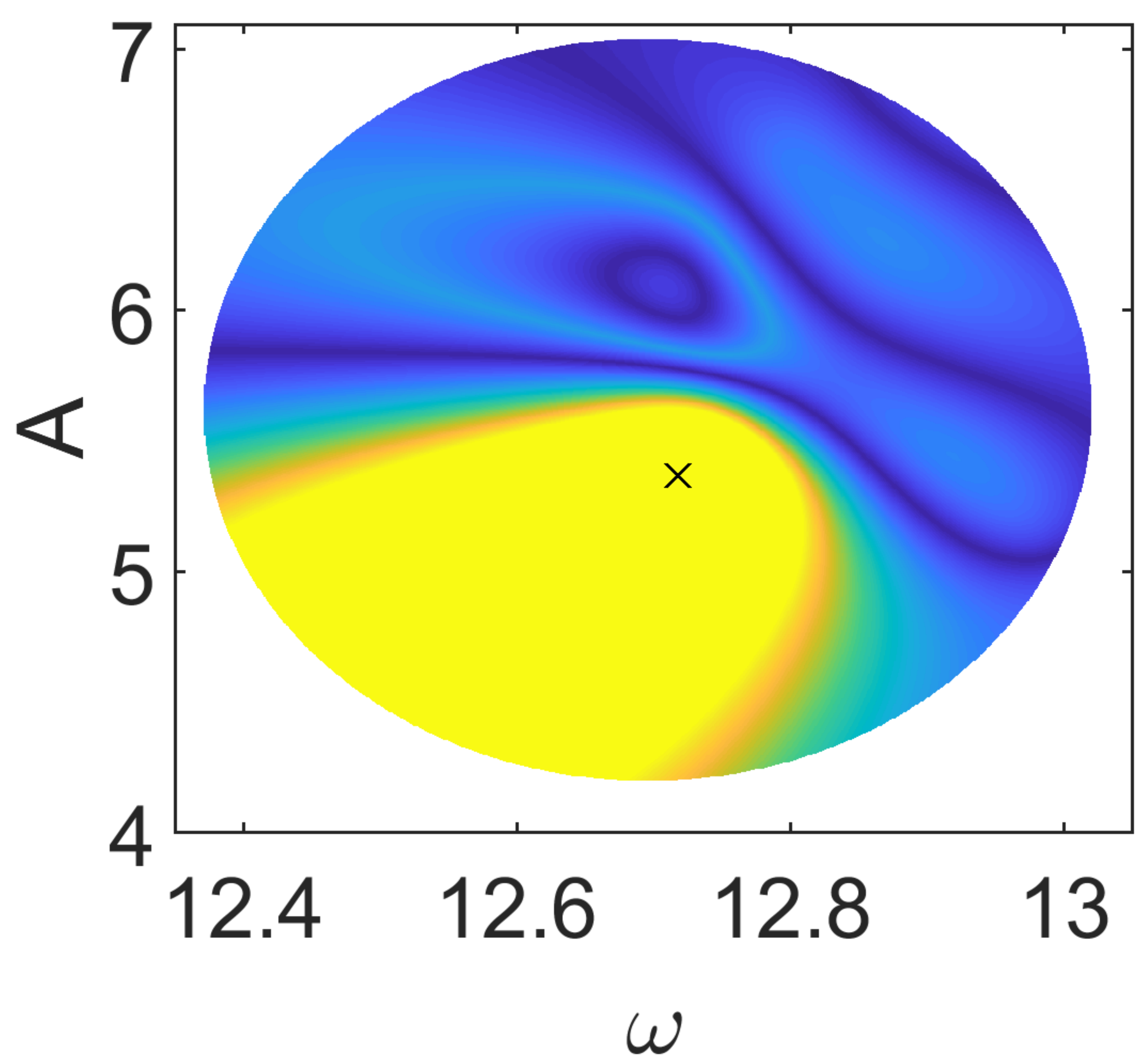}} &
\subfloat[]{\label{sens2}\includegraphics[width=0.32\textwidth]{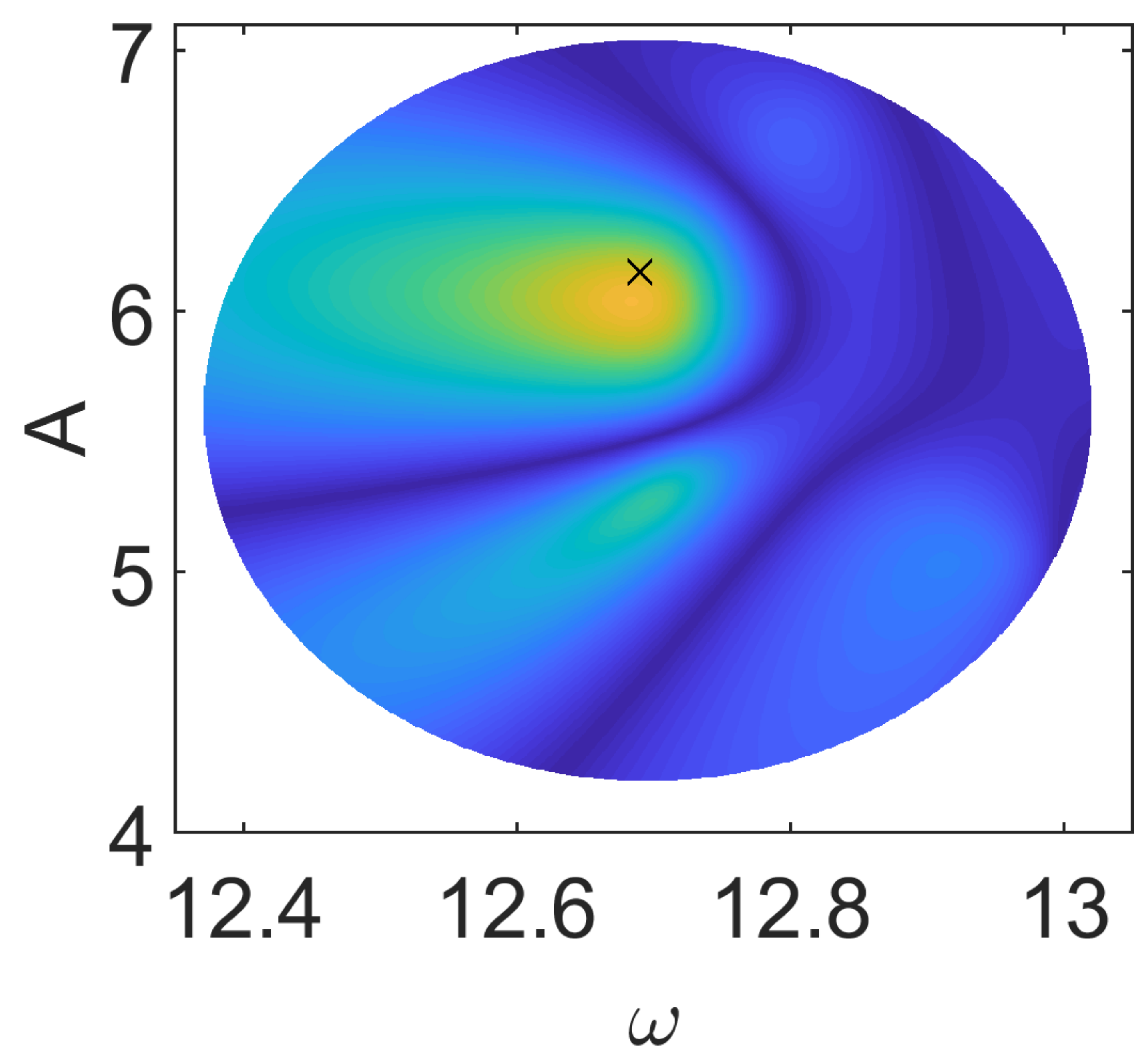}} &
\subfloat[]{\label{sens3}\includegraphics[width=0.32\textwidth]{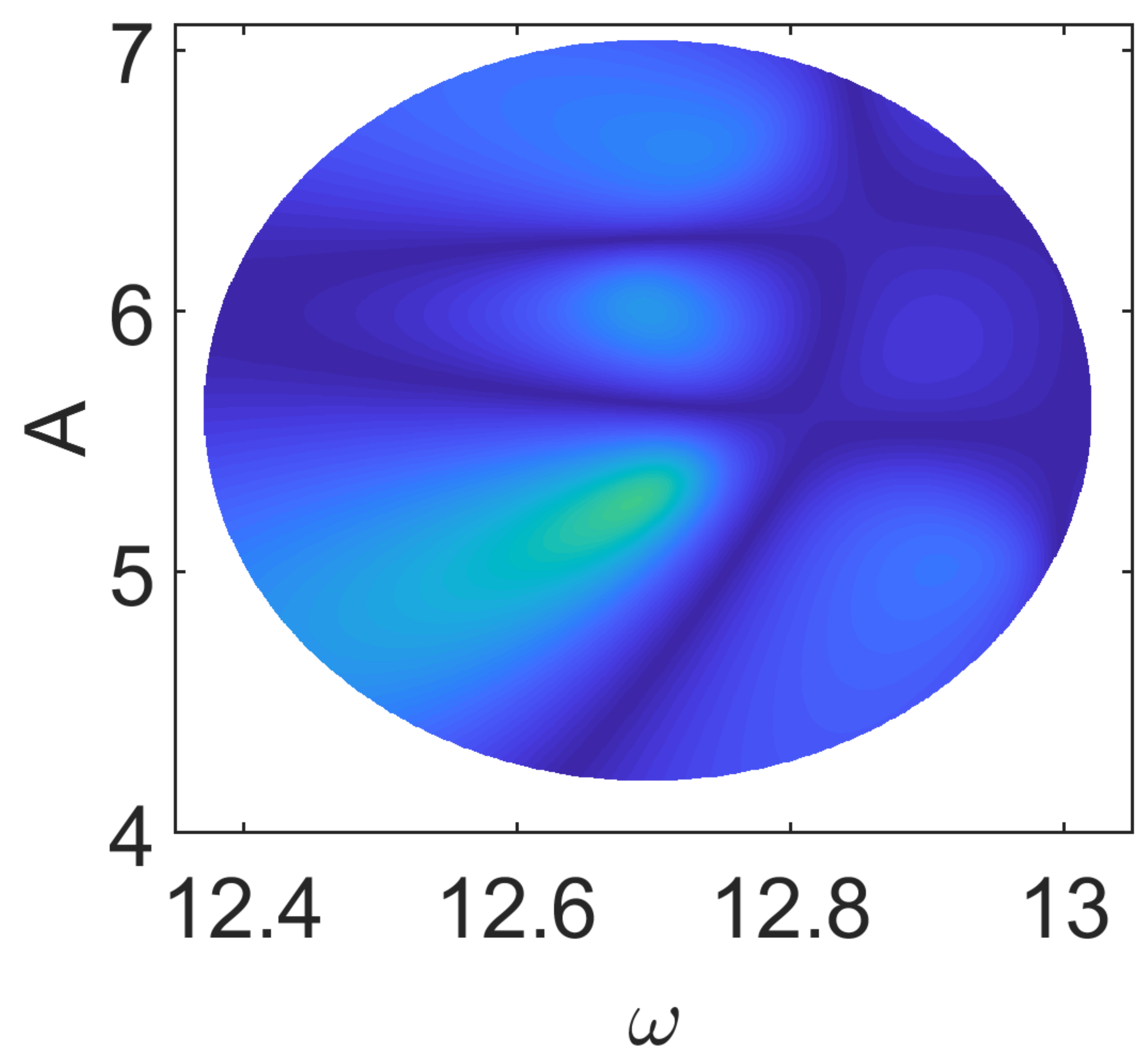}} 
\end{tabular*}
\begin{tikzpicture}
\begin{axis}[
    hide axis,
    scale only axis,
    height=0pt,
    width=0pt,
    colorbar horizontal,
    point meta min=0,
    point meta max=0.07,
    colorbar style={ylabel = {$\beta$}, width=10cm, xtick={0, 0.04, 0.07}}]
    \addplot [draw=none] coordinates {(0,0)};
\end{axis}
\end{tikzpicture}
\caption{Illustrates the new data selection process. The colourmap shows the sensitivity of the zero-problem $\beta$ as estimated by Eq.~\eqref{eq:comp}. (a -- c) Show the evolution of Eq.~\eqref{eq:comp} when new data points (\textcolor{black}{$\boldsymbol{\times}$}) are added to the GPR model. $4\times 10^{-2}$ is the user-defined threshold below which no additional data point is collected.}
\label{fig:point_selection}
\end{figure*}
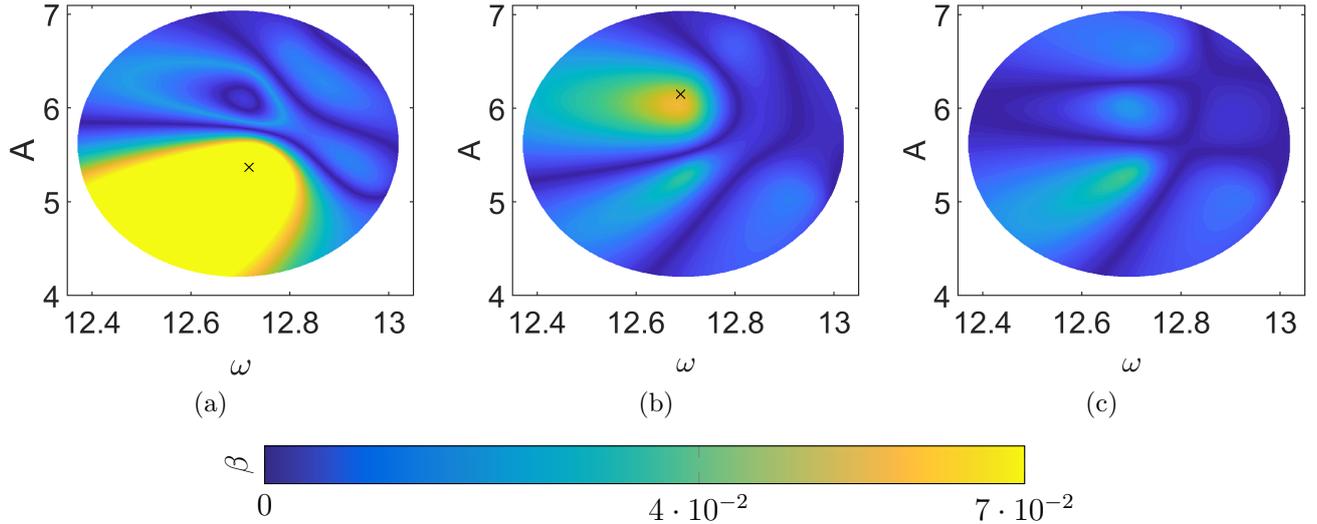

Figure~\ref{fig:point_selection} illustrates this data selection approach. The colourmap represents the sensitivity, $\beta$, of the zero-problem~\eqref{eq:LP_cond} to a new data point. Starting from Figure~\ref{fig:point_selection}(a), a new data point is collected where Eq.~\eqref{eq:comp} is the largest (\textcolor{black}{$\boldsymbol{\times}$}). This new data point is added to the GPR model and the correction step of the continuation algorithm repeated. Once a solution is found, the influence of a new data point on the new GPR model is again assessed (Figure~\ref{fig:point_selection}(b)). The region where the first additional data point was added is now observed to be significantly less influenced by any new point. According to the model predictions, informative data points are now located in another region where a second data point is eventually recorded. Following the same procedure, data points are collected until Eq.~\eqref{eq:comp} is below a user-defined tolerance across the whole region as in Figure~\ref{fig:point_selection}(c). At this stage, the solution of the continuation problem is said to be robust to new data and the continuation algorithm can perform a new prediction step. As the continuation algorithm progresses in parameter space, more data points are added to the GPR model. To keep computational costs low, an overall maximum number of data points $n_{\text{max}}$ in the model is maintained by removing data points that have less influence on the zero-problem.\\

During the experiment, the set of candidate data points is usually limited to 50 points. This is a much smaller set of points than the one considered for the colour maps in Figure~\ref{fig:point_selection}. This explains why in Figure~\ref{fig:point_selection}(b) there exist a small difference between the apparent location of the maximum of Eq.~\eqref{eq:comp} and the location where the new data point has been collected (\textcolor{black}{$\boldsymbol{\times}$}).\\

Other approaches to decide where to collect data points could have also been used. For instance, the effect of new data points on the solution of the zero-problem was investigated. However, this approach was found to give similar results to the method above while needing a solution to the nonlinear continuation problem for each candidate data point and hence being computationally much more expensive. Another approach was to select the data points for which the variance, $\text{var}[\Gamma_{\star}]$, of the predicted distribution was the largest. This approach was however discarded as new data points were being positioned at the periphery of the data set.\\ 

\section{Description of the experimental set-up}\label{sec:rig}
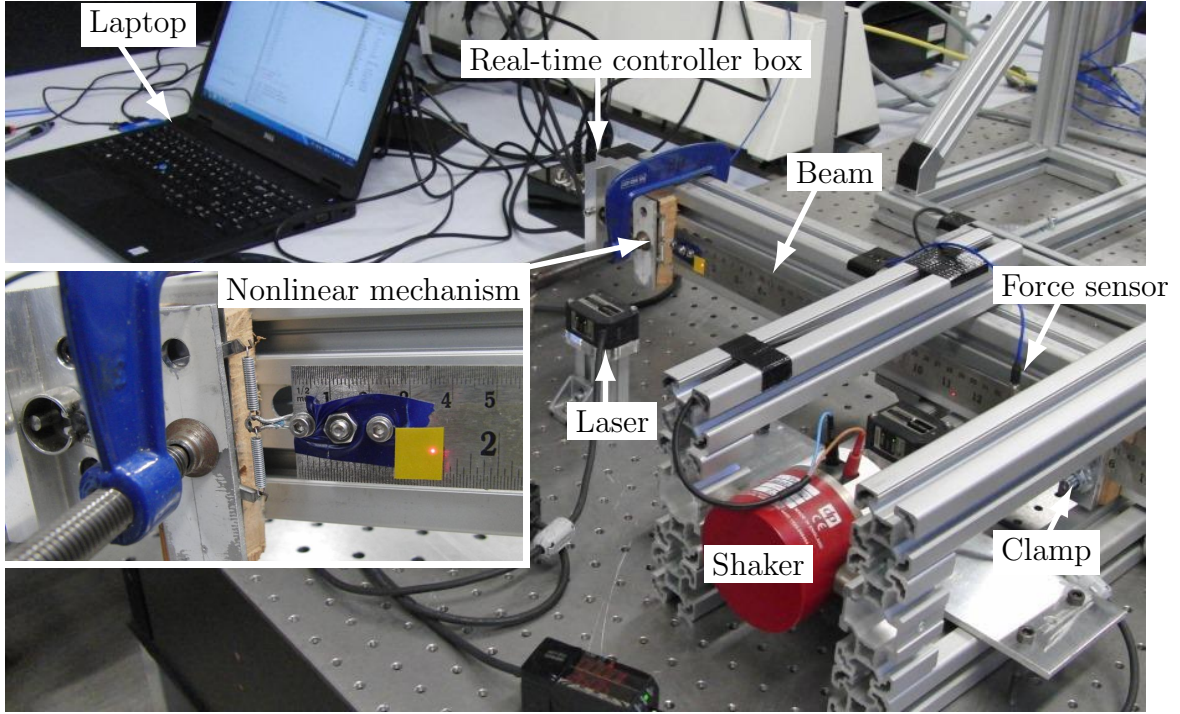
\begin{figure*}[t]
\centering
\begin{picture}(15, 9.5)
\fboxrule=0pt
\put(0, 0){\includegraphics[width=15cm]{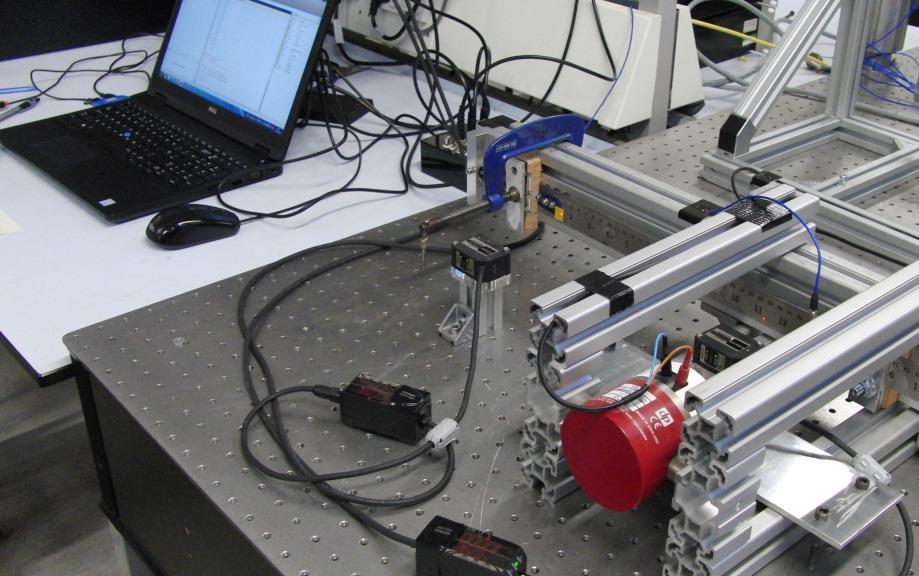}}
\put(-0.1, 2.0){\fcolorbox{white}{white}{\includegraphics[width=0.40\textwidth]{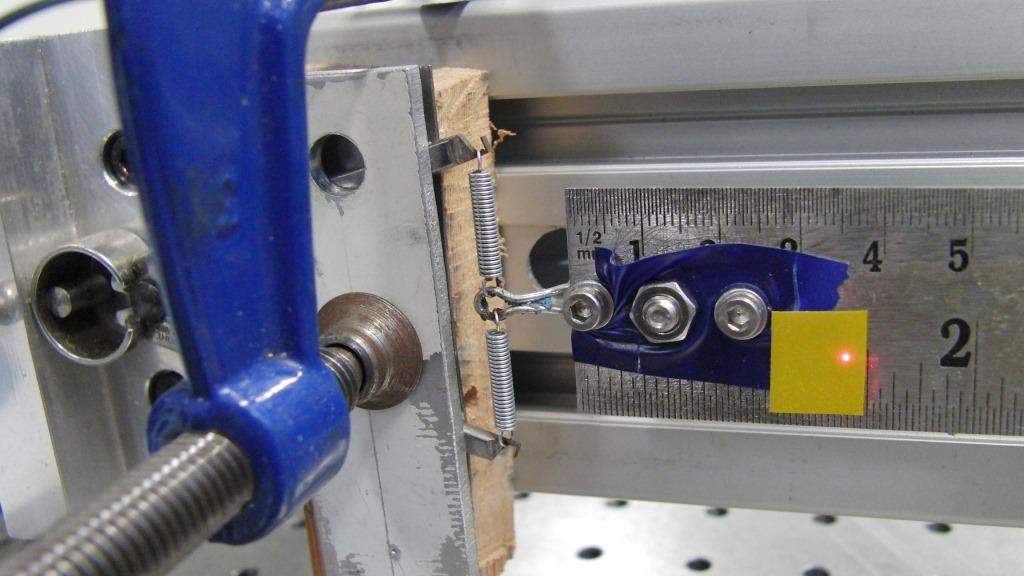}}}
\put(1.,9.0){\colorbox{white}{Laptop}}
\put(6.,8.5){\colorbox{white}{Real-time controller box}}
\put(2.8,5.45){\colorbox{white}{Nonlinear mechanism}}
\put(10.3,7.0){\colorbox{white}{Beam}}
\put(9.2,1.8){\colorbox{white}{Shaker}}
\put(7.4,3.7){\colorbox{white}{Laser}}
\put(13,5.5){\colorbox{white}{Force sensor}}
\put(13,2){\colorbox{white}{Clamp}}
\linethickness{2pt}
\put(1.5,8.85){\color{white}\vector(1,-1.5){.7}}
\put(7.8,8.5){\color{white}\vector(0,-1){1.1}}
\put(6.8,5.87){\color{white}\vector(7,1.5){1.7}}
\put(10.7,7.0){\color{white}\vector(-1,-2){0.6}}
\put(13.9,5.5){\color{white}\vector(-1,-2.){.5}}
\put(7.9,4){\color{white}\vector(0.,1.0){.8}}
\put(13.7,2.3){\color{white}\vector(1,2){.35}}
\end{picture}
\caption{Experimental set-up. The structure corresponds to a cantilever beam with a nonlinear mechanism attached at its free end. The structure is excited and controlled by means of an electro-dynamic shaker. A displacement laser sensor is used to measure the motion of the beam tip. CBC algorithms run on a laptop interconnected with the real-time controller box.}
\label{fig:setup}
\end{figure*}

The experimental set-up considered for the experimental demonstration of our algorithm is shown in Figure~\ref{fig:setup}. The main structure is a steel ($\rho \approx$ 7850 Kg/m$^3$, $E \approx$ 210 GPa) cantilever beam of length 380 mm, width 31.9 mm and thickness 1.9 mm. The beam's free end is attached to two linear springs arranged as shown in Figure~\ref{fig:setup}(b). This mechanism gives rise to geometric nonlinearity at large amplitudes. Previous work has shown that the stiffness properties of this mechanism can be approximated by a linear plus cubic term~\cite{Shaw16}. However, a mathematical model of the nonlinearity is unnecessary for CBC. As such, neither the identification of the nonlinear parameters nor the exploitation of the mathematical form of the nonlinearity were used. The length of the beam as well as the pre-tension in the springs were carefully adjusted such that the ratio between the natural frequencies of the first two bending modes is close to, but larger than 3. This leads to the presence of a 3:1 modal interaction between these modes.\\

The structure is excited approximately 40 mm away from the clamp using a GW-V4 Data Physics shaker powered by a PA30E amplifier. The force applied to the structure is measured using a PCB208C03 force transducer connected to a Kistler signal conditioner (type 5134). The vibrations of the beam are measured at the tip using an Omron ZX2-LD100 displacement laser sensor. The beam structure, the first laser sensor, the shaker and its power amplifier constitute the nonlinear experiment tested using CBC.\\

The algorithm used by the CBC method and presented in Section~\ref{sec:theory} is run on a laptop computer directly connected to the real-time controller (RTC) box via a USB cable. The RTC box consists of a BeagleBone Black on which the feedback controller used by CBC is implemented. Note that CBC algorithms do not run in real-time, only the feedback controller does. The BeagleBone Black is fitted with a custom data acquisition board (hardware schematics and associated software are open source and freely available~\cite{CBC_hardware}). All measurements are made at 1\,kHz sampling with no filtering. Estimations of the Fourier coefficients of the response, input force, and control action are calculated in real time on the control board using recursive estimators~\cite{Renson16}; however, this was for convenience rather than a necessity. \\

The $z$-domain transfer function of the controller used by the CBC technique is given by
\begin{equation}
\frac{U(z)}{E(z)} = \frac{0.0053}{z^3-2.4521 z^2+1.9725 z - 0.5155},
\label{eq:controller}
\end{equation} and aims to reduce the error $E(z)$ between the beam tip response, $y$ (laser 1), and a control reference signal, $y^*$, (see Section~\ref{sec:results_basic}). The control law, which was found to stabilise the dynamics of the experiment throughout the range of parameters considered in this study, was designed using pole-placement techniques and a linear model of the experiment. This model was obtained using low-amplitude broadband input-output data and captures the first two bending modes of the beam whose natural frequencies (damping ratios) were estimated at 11.49 Hz (0.026) and 36.45 Hz (0.022), respectively. Additional details on the derivation of the controller can be found in~\cite{Renson19} where the nonlinear frequency response curves of the present experiment were also investigated using CBC. Note that errors in the model do not affect the results as long as the model is sufficiently accurate for designing a stabilising feedback controller.\\

\section{Experimental results}\label{sec:results}
The CBC technique used in this paper to collect experimental data points is briefly reviewed in Section~\ref{sec:results_basic} and exploited to map out the complete response surface of the system of interest. The new online regression-based algorithm is then demonstrated in Section~\ref{sec:results_GPcont} where it is used to track fold points while the experiment is running. To validate the new algorithm, obtained results are directly compared with the fold curves that can be extracted by post-processing the response surface obtained using standard CBC approach as described in Section~\ref{sec:results_basic}.

\subsection{Extracting periodic responses}\label{sec:results_basic}
Considering periodic responses of the experiment, it is assumed that the control target, $y^*(t)$, and the control signal, $u(t)$, can be decomposed into a finite number $m$ of Fourier modes as
\begin{eqnarray}\label{eq:fourier_decomp}
    y^*(t) &=& \frac{A_0^*}{2} + \sum_{j=1}^{m} A_j^* \cos(j \omega t) + B_j^* \sin(j \omega t),\\
    u(t) &=& \underbrace{A^u_1 \cos (\omega t) + B^u_1 \sin (\omega t)}_{\text{fundamental}} \\
    & & + \underbrace{\frac{A_0^u}{2} + \sum^{m}_{j=2} A_j^u \cos (j \omega t) + B_j^u \sin (j \omega t)}_{\text{higher harmonics}}.
    \label{eq:FourierExp}
\end{eqnarray}
CBC seeks a reference signal $y^*(t)$ for which the control signal $u$ is equal to zero for all times. When this is achieved, the control signal is non invasive and does not alter the position in parameter space of periodic responses compared to the underlying uncontrolled system of interest. To achieve this non-invasive control signal, the higher harmonics $(A^u_0, A_j^u, B_j^u)_{j=2}^m$ have to be eliminated by finding suitable reference signal coefficients $(A^*_0, A_j^*, B_j^*)_{j=2}^m$. This can be performed using a derivative-free Picard-iteration algorithm. As for the fundamental harmonic component of the control signal, $A^u_1 \cos (\omega t) + B^u_1 \sin (\omega t)$, it can be viewed as the harmonic excitation applied to the experiment. If no other external excitation is applied to the experiment, the total harmonic excitation applied to the experiment is simply given by $\Gamma = \sqrt{A_1^{u^2}+B_1^{u^2}}$. The reader is referred to~\cite{Barton13,Renson17} for a more detailed discussion on this approach that is often referred to as the simplified CBC method.\\

The higher-harmonic coefficients $(A^*_0, A_j^*, B_j^*)_{j=2}^m$ are determined to cancel out the higher-harmonic coefficients of the control signal. One of the fundamental target coefficients (here $B_1^*$) can also be set equal to zero in order to set a phase reference between the response and the excitation. As such, the only two adjustable inputs to the experiment are the frequency of excitation $\omega$ and the fundamental target coefficients $A_1^*$. The amplitude of the excitation $\Gamma$ is viewed as a free parameter as it is not fully determined by the user but depends on the response $y(t)$ and target $y^*(t)$.\\

\begin{figure*}[t]
\centering
\includegraphics[width=0.9\textwidth]{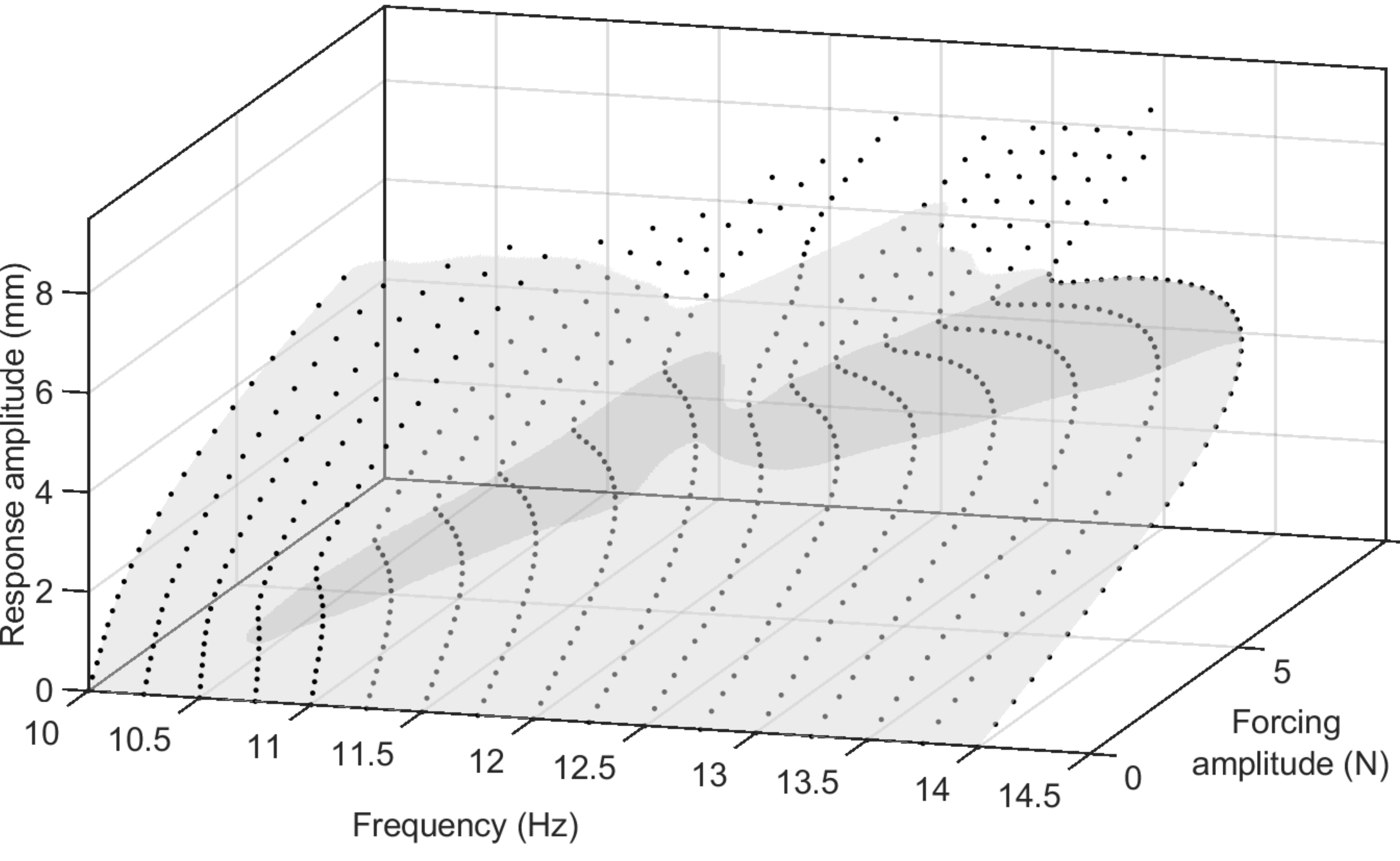}
\caption{$S$-curves collected using the simplified CBC algorithm presented in Section~\ref{sec:results_basic}.  Gray surface is the output of a GPR model including all recorded data points (\textcolor{black}{$\boldsymbol{\bullet}$}). The darker gray region is delimited by a fold curve. Data points in this region correspond to unstable response of the underlying uncontrolled experiment.}
\label{fig:RefSurf}
\end{figure*}

Considering a constant frequency, the response of the experiment to harmonic input can be investigated by increasing the free fundamental target coefficient $A_1^*$. The response curves collected in this way form curves that resemble a $S$, as shown in Figure~\ref{fig:RefSurf}. At each data point (\textcolor{black}{$\boldsymbol{\bullet}$}), full time series measurements containing 2000 samples, or between 22 and 29 oscillation cycles, were made and used to estimate response and force amplitudes. The number of Fourier modes considered in the calculation of amplitudes is equal to seven throughout the rest of the paper.\\

Response curves in Figure~\ref{fig:RefSurf} were collected every 0.25 Hz between 10 Hz and 14 Hz with 0.2 mm steps in $A_1^*$. The $S$-curve obtained at 12.5 Hz is notably different from the other $S$-curves and includes an additional inflection point in its upper part. This curve lies in the region where a modal interaction occurs between the first two modes of the structure and marks a change in the form of the response surface of the system. This change in the response surface leads to the presence of an isola as further discussed in Section~\ref{sec:results_GPcont}.\\

The $S$-curves collected can be post-processed to extract fold curves which will then be compared to the fold curves directly tracked during the experiment. The post-processing approach used to extract these curves is similar in principle to the one used in the online algorithm but with the clear distinction that it does not influence the data points collected. A GPR model including all collected data points is created and numerical continuation is exploited to follow the solutions of Eq.~\eqref{eq:LP_cond}. Figure~\ref{fig:RefSurf} shows in gray the continuous surface constructed from the GPR model. The darker gray region indicates where the underlying uncontrolled experiment is unstable and thus unobservable without stabilising feedback control. The boundary of that region is the fold curve found using continuation. This curve is sensitive to the spacing between the different S-curves as well as the choice of hyper-parameters for the GPR. In particular, the fold curve presents artifact oscillation as discussed later. The online algorithm will overcome these issues.\\

\subsection{Demonstration of the online regression-based algorithm}\label{sec:results_GPcont}
\begin{figure*}[htp]
\centering
\begin{tabular*}{1.\textwidth}{@{\extracolsep{\fill}} c c}
\multicolumn{2}{c}{\subfloat[]{\label{start1}\includegraphics[width=0.9\textwidth]{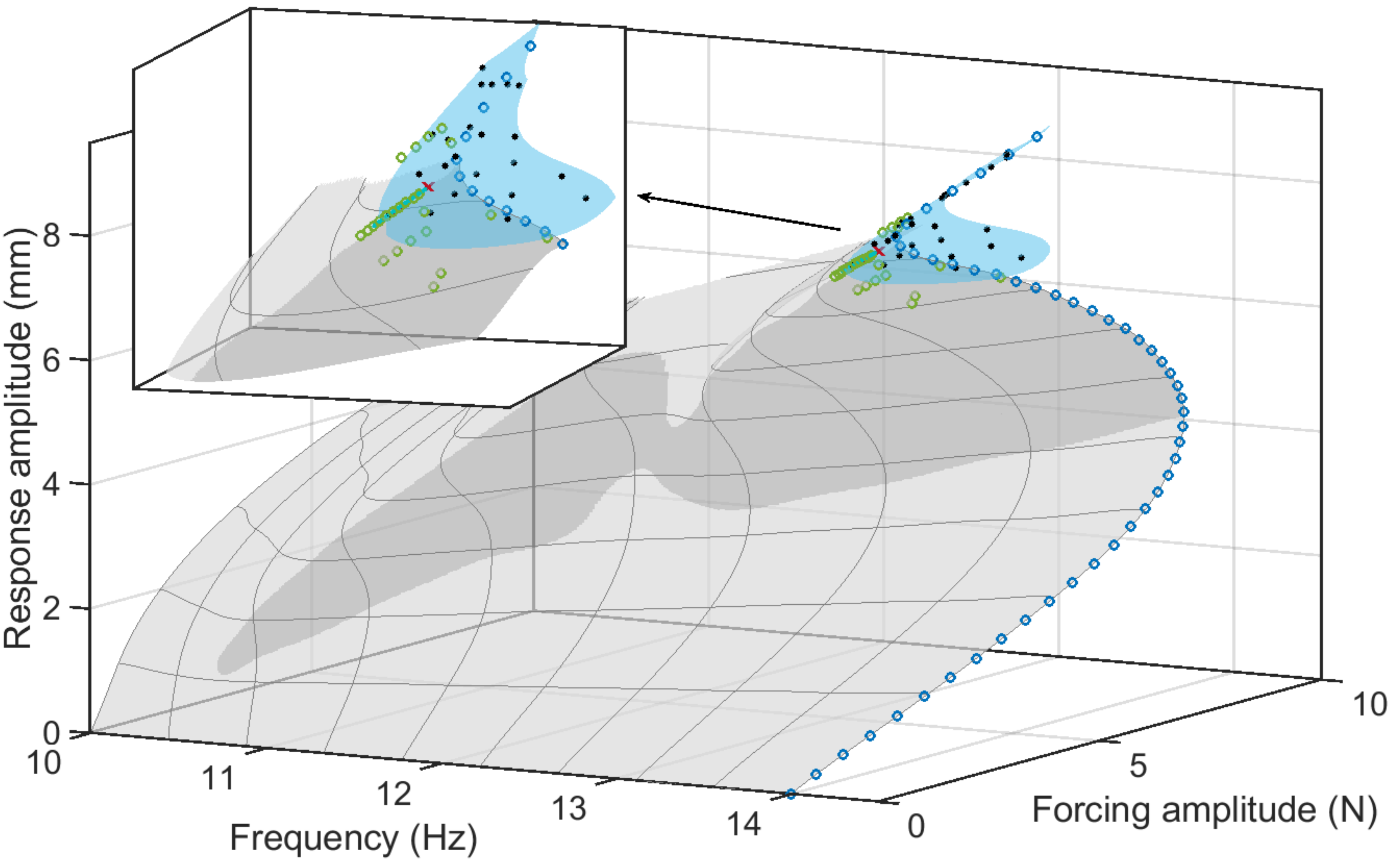}}}\\
\subfloat[]{\label{5stepsa}\includegraphics[width=0.45\textwidth]{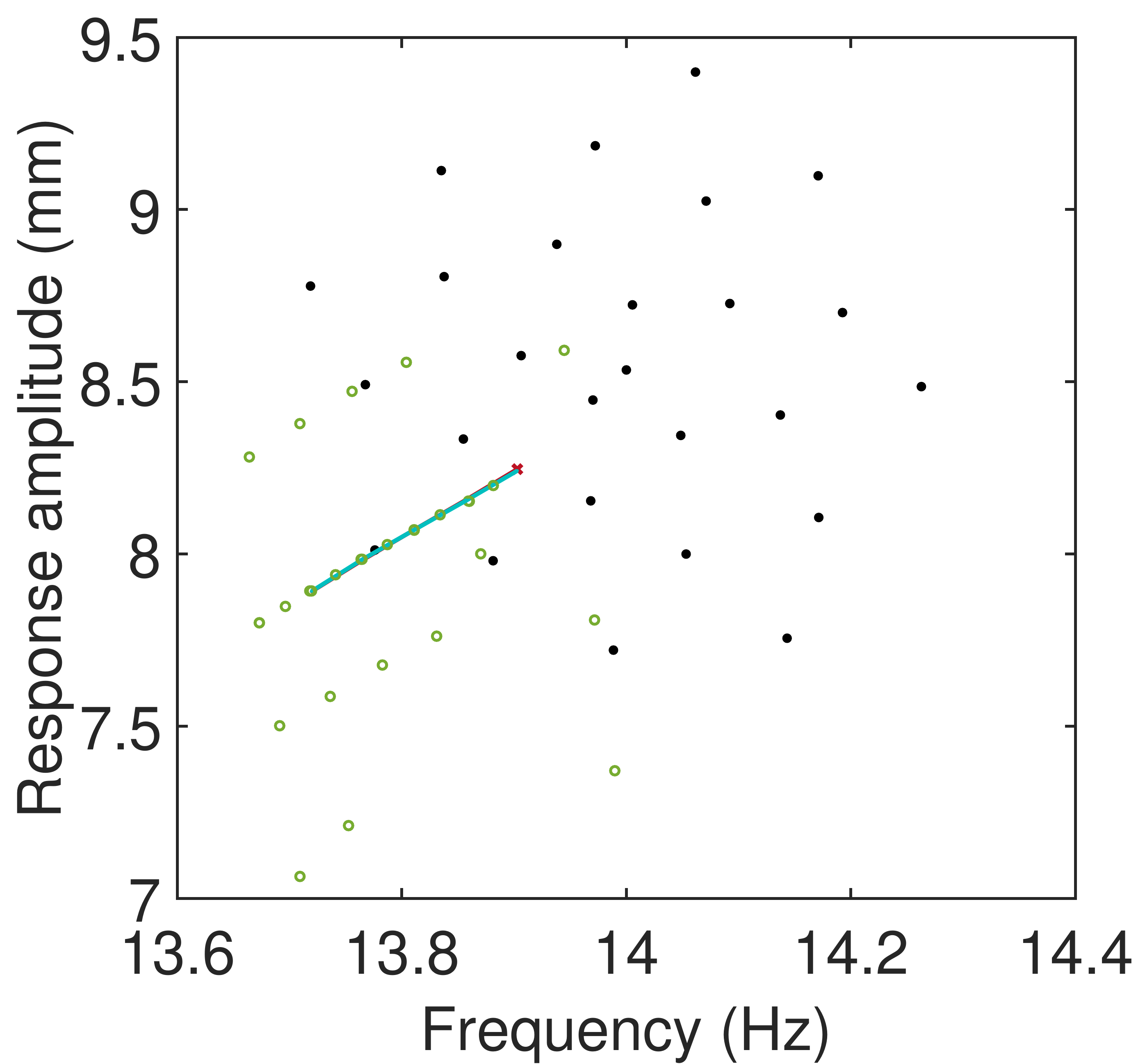}} &
\subfloat[]{\label{5stepsb}\includegraphics[width=0.45\textwidth]{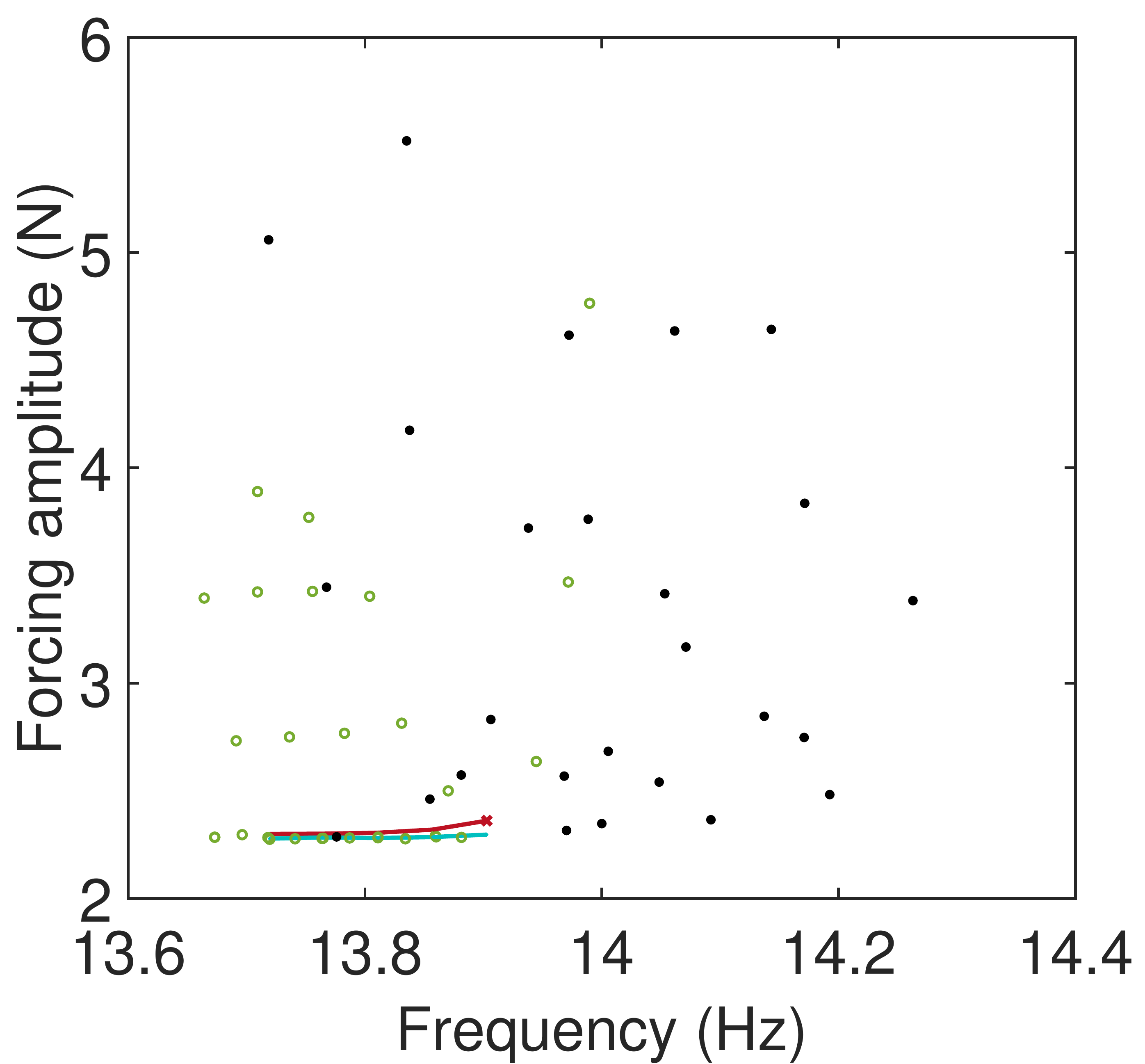}}
\end{tabular*}
\caption{First four steps of the regression-based continuation algorithm. (\textcolor{dblue}{$\boldsymbol{\circ}$}) S-curve used to initialize the algorithm close to a fold. Data points (\textcolor{black}{$\boldsymbol{\bullet}$}) used to estimate hyper-parameters and find the first fold point. Data points (\textcolor{dgreen}{$\boldsymbol{\circ}$}) collected during the continuation process. (\textcolor{dred}{$\boldsymbol{-\bullet-}$}) and (\textcolor{dgreen}{$\boldsymbol{-\bullet-}$}) are the fold curves found using continuation and where the force amplitude is estimated by the GPR model and measured experimentally, respectively. (\textcolor{lblue2}{$\blacksquare$}) First local GPR model. To help visualization, the response surface obtained in Figure~\ref{fig:RefSurf} is superimposed to the data obtained with the new algorithm. (b, c) Show two-dimensional projections of the data points and curves shown in (a).}
\label{fig:first_steps}
\end{figure*}

The online regression-based continuation algorithm presented in Section~\ref{sec:theory} is now demonstrated. Figure~\ref{fig:first_steps} illustrates the first four steps of the algorithm which is initialized in the neighbourhood of a fold found at high response amplitudes on one of the $S$-curves collected in Section~\ref{sec:results_basic} (\textcolor{dblue}{$\boldsymbol{\circ}$}). Around that starting point, 25 new data points (\textcolor{black}{$\boldsymbol{\bullet}$}) regularly distributed in frequency and control target amplitude are collected. This first data set is then used to estimate the hyper-parameters of the GPR model using the marginalization approach presented in Section~\ref{sec:algo_gpr}. Hyper-parameters were found to be $\boldsymbol{\theta} = (\sigma_n^2, \sigma_f^2, l_{\omega}, l_A) \approx (0.02, 2.02, 0.30, 1.09)$. The noise variance, $\sigma^2_n$, which is found to be relatively small compared to the other hyper-parameters, captures the effects of measurement noise. It also encompasses potential inaccuracies that could arise from the inaccurate cancellation of the control signal higher harmonics during data collection. \\

\begin{figure}[t]
\centering
\includegraphics[width=0.45\textwidth]{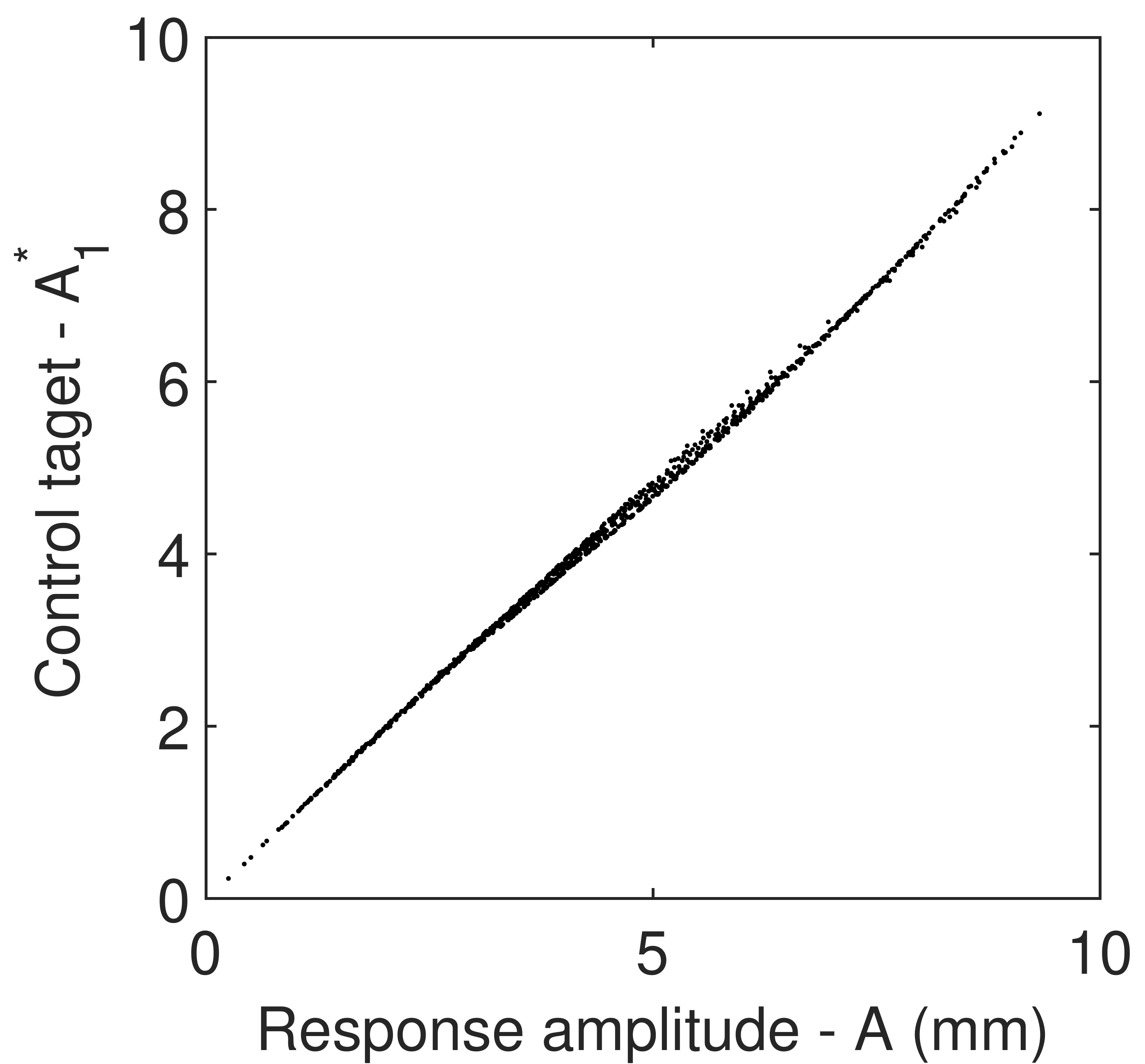}
\caption{Mapping between the response amplitude $A$ (continuation variable) and the fundamental target coefficient $A_1^*$ (experimental input). All the data points (\textcolor{black}{$\boldsymbol{\bullet}$}) presented in this paper have been included.}
\label{fig:mapping}
\end{figure}

After estimation of the hyper-parameters, the first GPR model is created (\textcolor{lblue2}{$\blacksquare$}) and a first solution to Eq.~\eqref{eq:LP_cond} is found (\textcolor{dred}{$\boldsymbol{\times}$}). The sequence of prediction and correction steps of the numerical continuation algorithm can then be started. Data points collected during the continuation process are shown in green (\textcolor{dgreen}{$\boldsymbol{\circ}$}) and appear to be regularly distributed around the fold of the response surface (see two-dimensional projections in Figure~\ref{fig:first_steps}(b, c)). This comes from the regular pattern used for the generation of candidate data points and the sensitivity criterion~\eqref{eq:comp} that is found to almost always have a structure similar to the one observed in Figure~\ref{fig:point_selection}(a, b). Note that to collect data points at the particular response amplitudes specified by the data selection method in Section~\ref{sec:algo_improve}, a second regression mapping $A$ to $A_1^*$ is created during the continuation. This mapping is observed to be almost linear throughout the range of parameters considered in this study (Figure~\ref{fig:mapping}).\\

In Figure~\ref{fig:first_steps}, the red curve (\textcolor{dred}{$\boldsymbol{- \bullet -}$}) represents the solutions found by the continuation algorithm and the force amplitude predicted at the solution by the GPR model. For each solution point, a data point is also collected (\textcolor{dgreen}{$\boldsymbol{- \bullet -}$}) and added to the GPR model. At the beginning of the algorithm (around 14 Hz), the forcing amplitude estimated by the GPR model at the solution is slightly different from the force amplitude actually measured. However, the GPR model predictions improve as more data points are added to it and the difference between the two curves (\textcolor{dred}{$\boldsymbol{- \bullet -}$}, \textcolor{dgreen}{$\boldsymbol{- \bullet -}$}) becomes increasingly smaller, and negligible after 3 continuation steps.\\

\begin{figure*}[tp]
\centering
\begin{tabular*}{1.\textwidth}{@{\extracolsep{\fill}} c c}
\multicolumn{2}{c}{\subfloat[]{\label{3D}\includegraphics[width=0.9\textwidth]{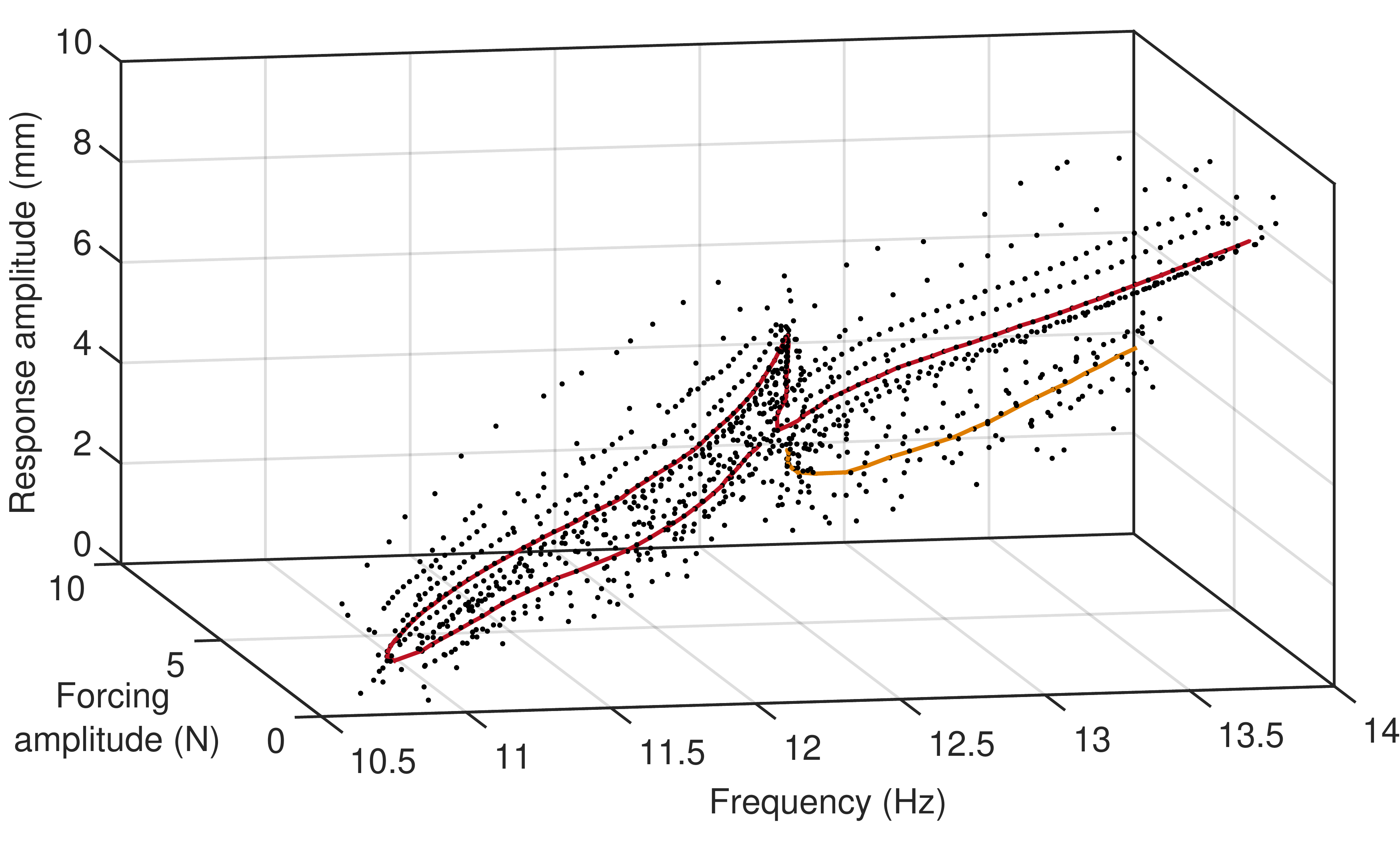}}}\\
\subfloat[]{\label{proj1}\includegraphics[width=0.45\textwidth]{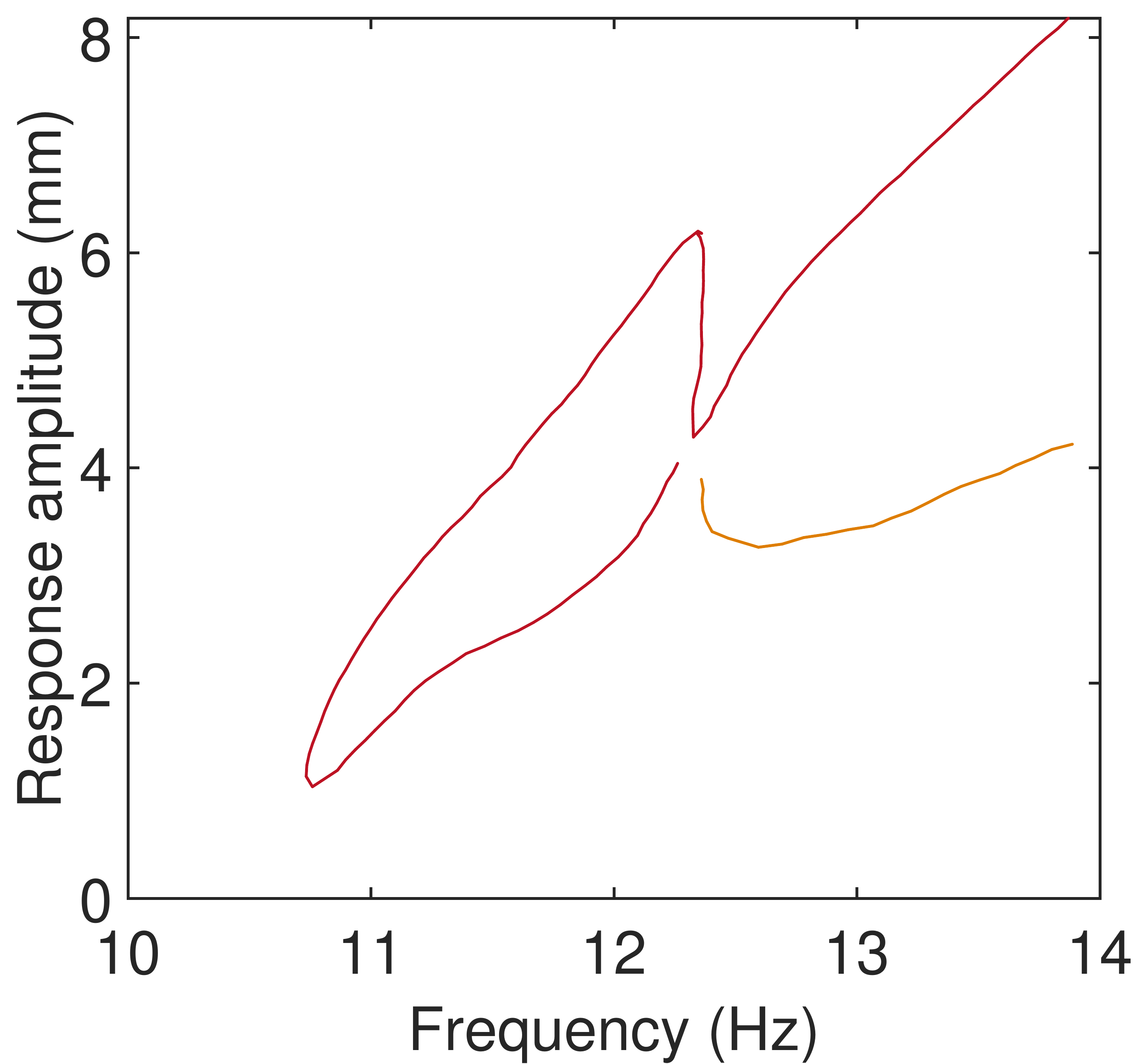}} &
\subfloat[]{\label{proj2}\includegraphics[width=0.45\textwidth]{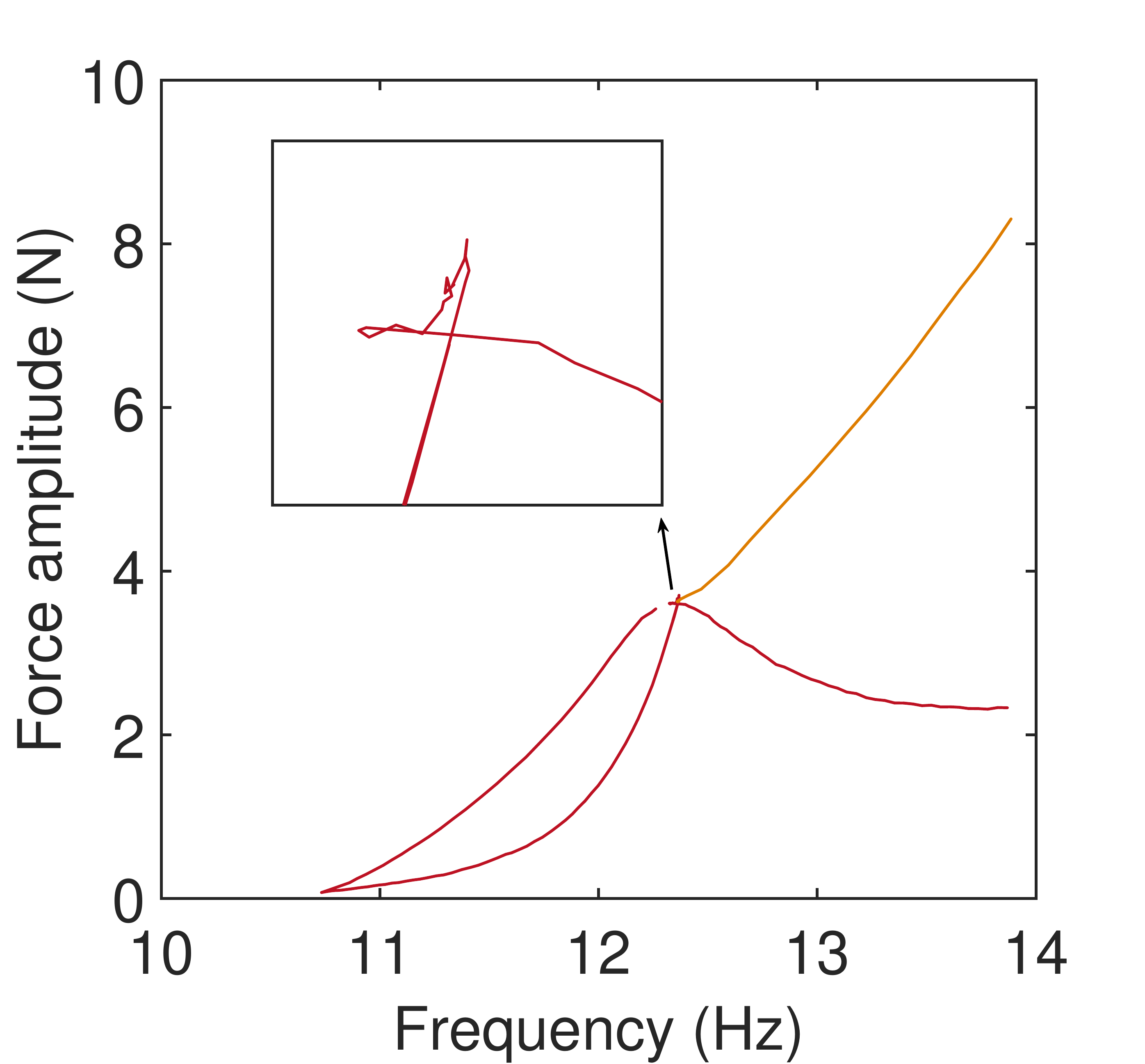}}
\end{tabular*}
\caption{(\textcolor{dred}{\textbf{---}}, \textcolor{orange}{\textbf{---}}) Fold curves captured experimentally using the regression-based continuation algorithm. (\textcolor{black}{$\boldsymbol{\bullet}$}) Data points automatically collected during continuation using the method proposed in Section~\ref{sec:algo_improve}. (b, c) Two-dimensional projections of the fold curves shown in (a). The close-up in (c) shows the presence of a swallowtail catastrophe.}
\label{fig:3D_LP_curve}
\end{figure*}

Figure~\ref{fig:3D_LP_curve} illustrates the results obtained when applying the continuation algorithm throughout the parameter range of interest. All the data points collected during the continuation are shown in black (\textcolor{black}{$\boldsymbol{\bullet}$}), and the curve of fold points obtained by following the curve previously shown in Figure~\ref{fig:first_steps} is in red (\textcolor{dred}{$\boldsymbol{-}$}). Two-dimensional projections of the fold curve are shown in Figure~\ref{fig:3D_LP_curve}(b, c). The curve has a complex geometry; especially in the region between 12 -- 13 Hz where it is observed to increase in response amplitude at almost constant forcing frequencies (Figure~\ref{fig:3D_LP_curve}(b)). This feature corresponds to the `saturation' of the resonance frequency commonly observed in nonlinear frequency response curves in the vicinity of modal interactions~\cite{Chen17b}. Closer inspection at this region reveals that the fold curve forms the shape of a swallowtail catastrophe (Figure~\ref{fig:3D_LP_curve}(c))~\cite{PostonBook}. Although this feature is small and appears to be more affected by measurement uncertainties than the rest of the curve, it was consistently observed in this region of the parameter space. After this interaction region, the algorithm was found to reach and successfully pass through a cusp (around 11 Hz).\\

The use of multivariate regression combined with a predictor-corrector continuation technique makes the online regression-based algorithm more general and robust than the bifurcation tracking algorithm proposed in~\cite{Renson17}. This latter standard method could not have tracked a curve that exhibits the complex features of the fold curve in Figure~\ref{fig:3D_LP_curve}. The present algorithm was however not able to systematically continue the fold curve across the modal interaction region. In particular, when the interaction region was approached from below, the fold curve was found to come very close to a previously measured section of the curve. Within uncertainty, this close proximity between the two portions of the curve gives rise to a branch point singularity which cannot be addressed by our current continuation algorithm. As such, the online algorithm was often found to `jump' to and continue the upper part of the fold curve. To capture the low-amplitude section of the fold curve located between 12.5 Hz and 14 Hz (\textcolor{dorange}{$\boldsymbol{-}$} in Figure~\ref{fig:3D_LP_curve}), the online algorithm was restarted around the second fold point observed at low amplitude on the $S$-curve in Figure~\ref{fig:first_steps}.\\

\begin{figure*}[t]
\centering
\begin{tabular*}{1.\textwidth}{@{\extracolsep{\fill}} c c}
\subfloat[]{\label{nlfr1}\includegraphics[width=0.45\textwidth]{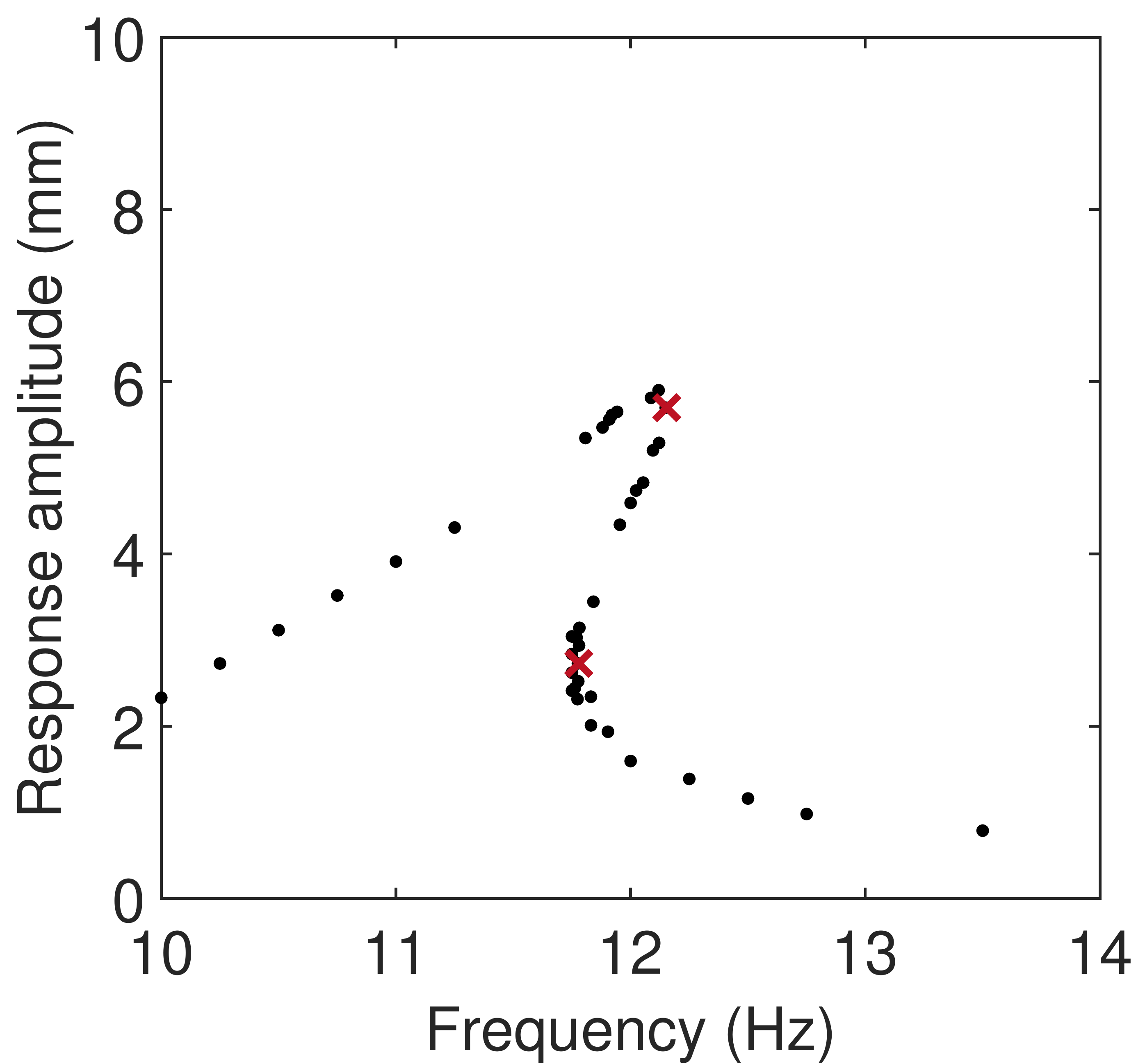}} &
\subfloat[]{\label{nlfr2}\includegraphics[width=0.45\textwidth]{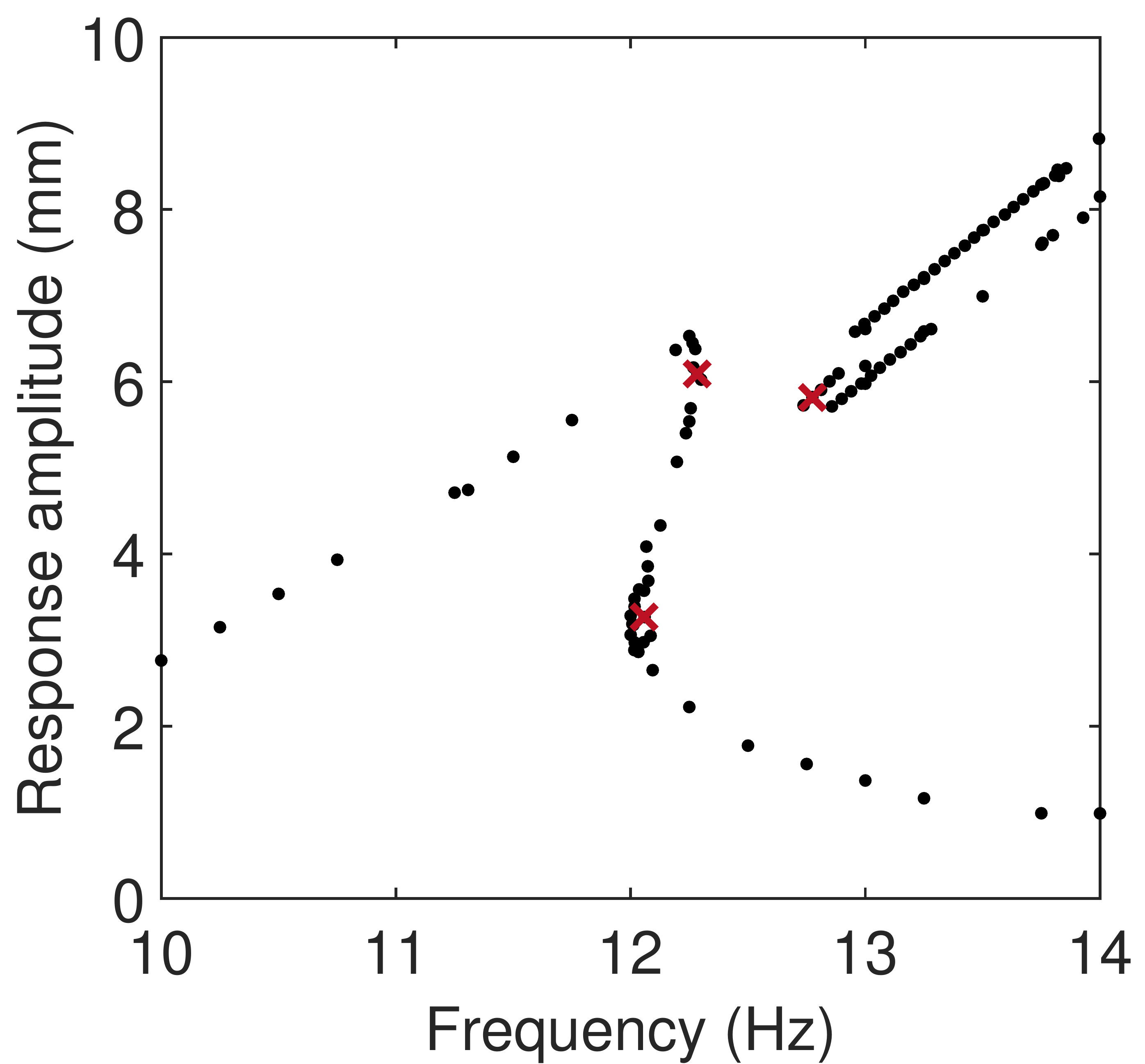}}
\end{tabular*}
\caption{Nonlinear frequency response of the system obtained by selecting data points in Figures~\ref{fig:RefSurf} and~\ref{fig:3D_LP_curve} that correspond approximately (within $\pm$5\%) to (a) 2.0 N and (b) 2.9 N. Fold points (\textcolor{dred}{$\boldsymbol{\times}$}) obtained with the regression-based continuation algorithm of Section~\ref{sec:theory}.}
\label{fig:NLFR}
\end{figure*}

\begin{figure*}[h]
\centering
\begin{tabular*}{1.\textwidth}{@{\extracolsep{\fill}} c c}
\subfloat[]{\label{comp1}\includegraphics[width=0.45\textwidth]{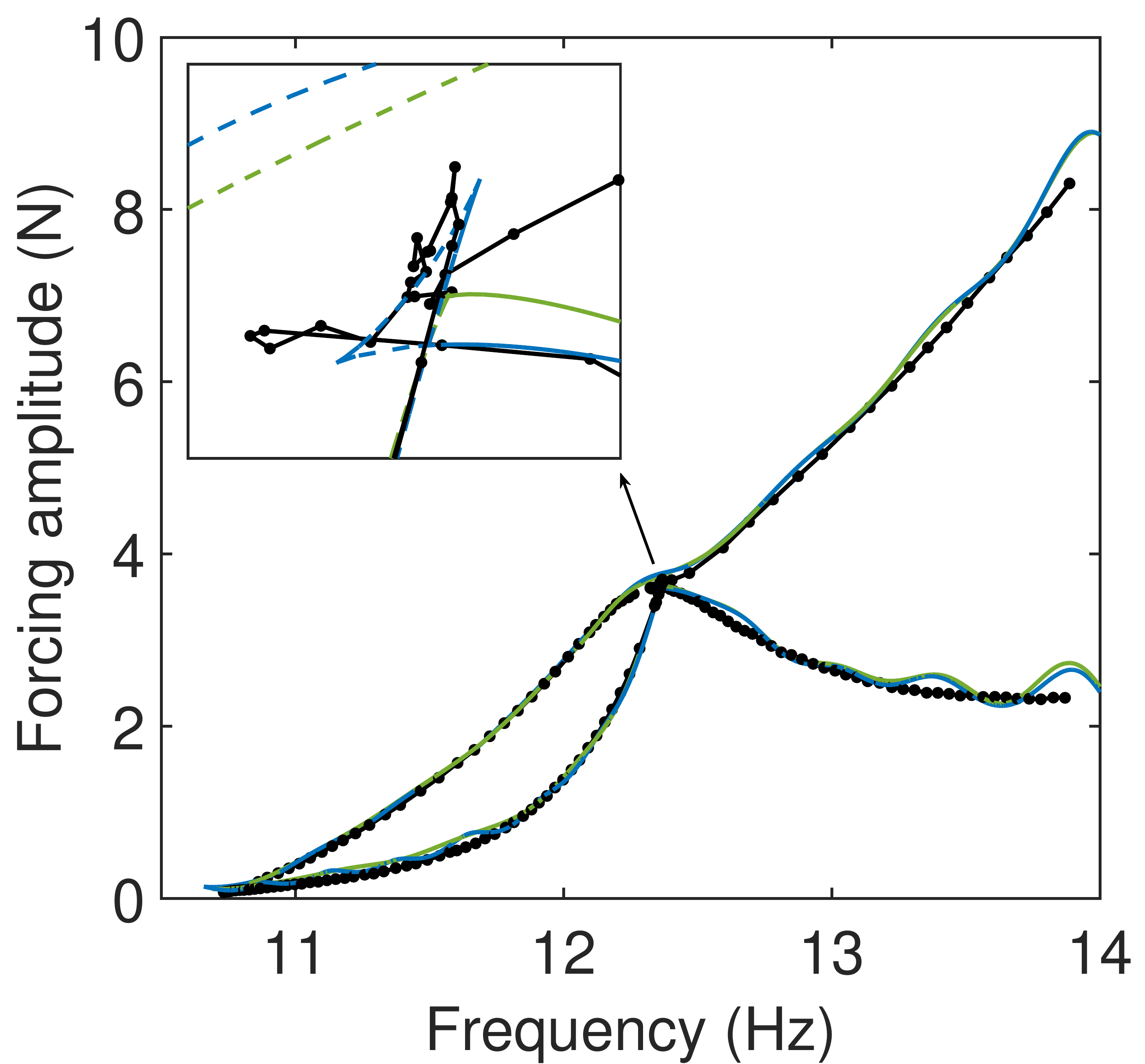}} &
\subfloat[]{\label{comp2}\includegraphics[width=0.45\textwidth]{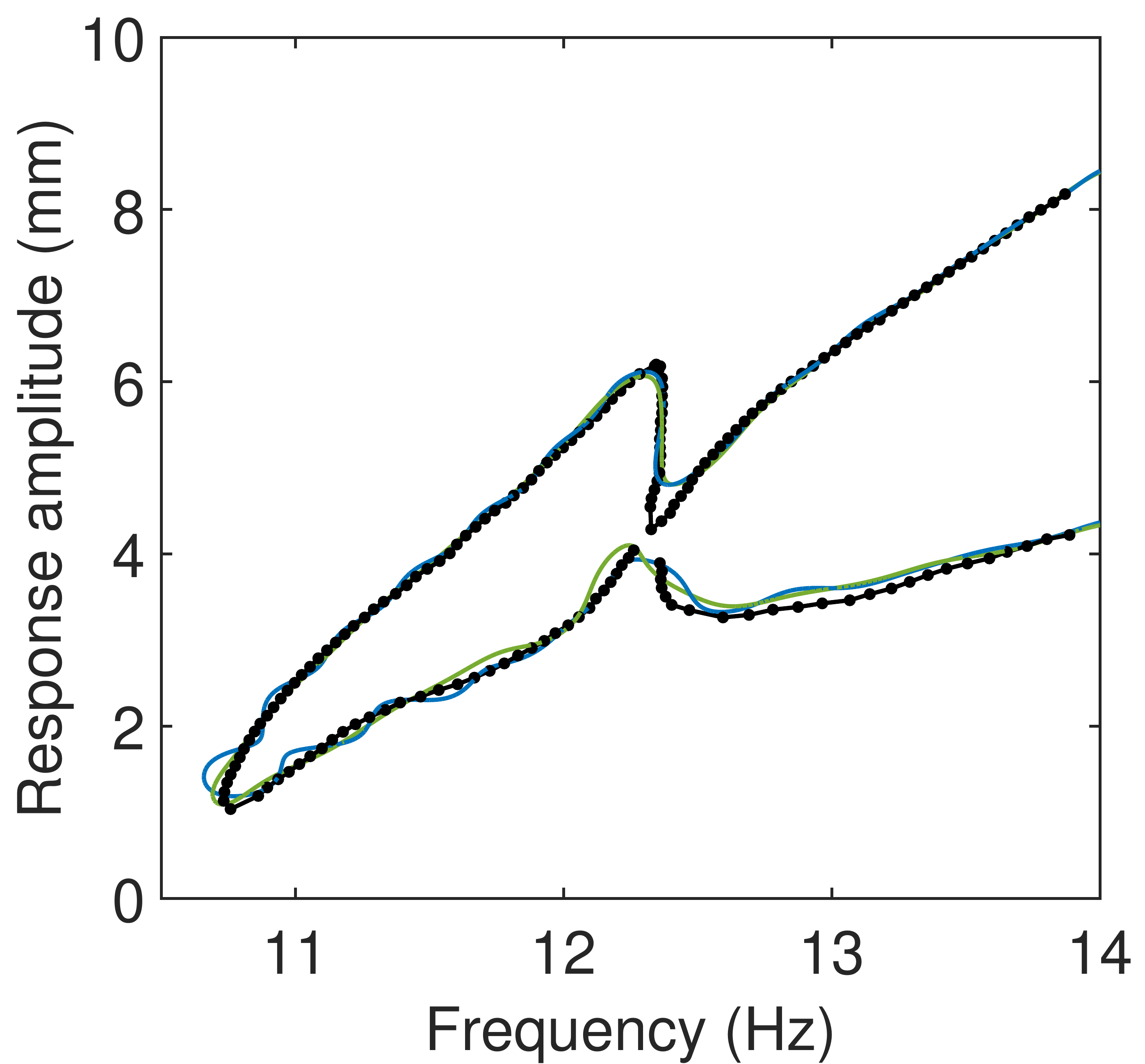}}
\end{tabular*}
\caption{Comparison between the fold curve (\textcolor{black}{$\boldsymbol{- \bullet -}$}) obtained during the experiment and the fold curves (\textcolor{dgreen}{$\boldsymbol{-}$}, \textcolor{dblue}{$\boldsymbol{-}$}) obtained by post-processing the data points of the response surface. Hyper-parameters used to obtain (\textcolor{dgreen}{$\boldsymbol{-}$}) are identical to those used in the experiment, whereas the hyper-parameters used for (\textcolor{dblue}{$\boldsymbol{-}$}) result from the maximisation of the marginal likelihood for all the data points in the response surface.}
\label{fig:vs_ref}
\end{figure*}

To validate the continuation results, fold points found at a particular force amplitude are compared to the nonlinear frequency response (NLFR) curve obtained for the same excitation level. To generate the NLFR curve, the raw data points shown in Figures~\ref{fig:RefSurf} and~\ref{fig:3D_LP_curve} that correspond approximately to the chosen force amplitude (within $\pm$5\%) are represented on the same graph. The NLFR curves obtained for 2.0 N and 2.9N are shown in Figures~\ref{fig:NLFR}(a) and~\ref{fig:NLFR}(b), respectively. Fold points are in red (\textcolor{dred}{$\boldsymbol{\times}$}). At 2N (Figure~\ref{fig:NLFR}(a)), the resonance peak leans towards higher frequencies which highlights the hardening character of the nonlinearity in the system. The locations of the two fold points found at that level match the folds in the NLFR. Data points located between them are observed to be unstable periodic responses of the underlying uncontrolled experiment. As the excitation is increased to 2.9 N (Figure~\ref{fig:NLFR}(b)), the resonance peak is qualitatively unchanged but an isolated branch of high-amplitude responses disconnected from the main resonance peak has appeared~\cite{Shaw16,Renson19}. The presence of this isola could have been inferred from Figure~\ref{fig:3D_LP_curve}(c) as a third fold point appears when the force amplitude becomes larger than approximately 2.1 N.\\

To further validate the fold curves obtained using the online regression-based algorithm, results are compared to the fold curve obtained by post-processing the data points of the full response surface. In particular, a GPR model including only the data points of Figure~\ref{fig:RefSurf} is created, and numerical continuation is applied to find and track the solutions of~\eqref{eq:LP_cond}. Two sets of hyper-parameters are considered to post-process the data. The first one corresponds to the set of hyper-parameters used by the online algorithm. The second one is obtained by maximising the marginal likelihood for all the data points of the response surface. Hyper-parameters found in this latter case are $(\sigma_n^2, \sigma_f^2, l_{\omega}, l_A) \approx (0.01, \;2.66,\allowbreak \;0.28, \;0.73)$. Although different from the hyper-parameters found at the initialization of the online continuation algorithm, they are found to have comparable magnitudes.\\

\begin{figure*}[h]
\centering
\begin{tabular*}{1.\textwidth}{@{\extracolsep{\fill}} c c}
\subfloat[]{\label{uq1}\includegraphics[width=0.45\textwidth]{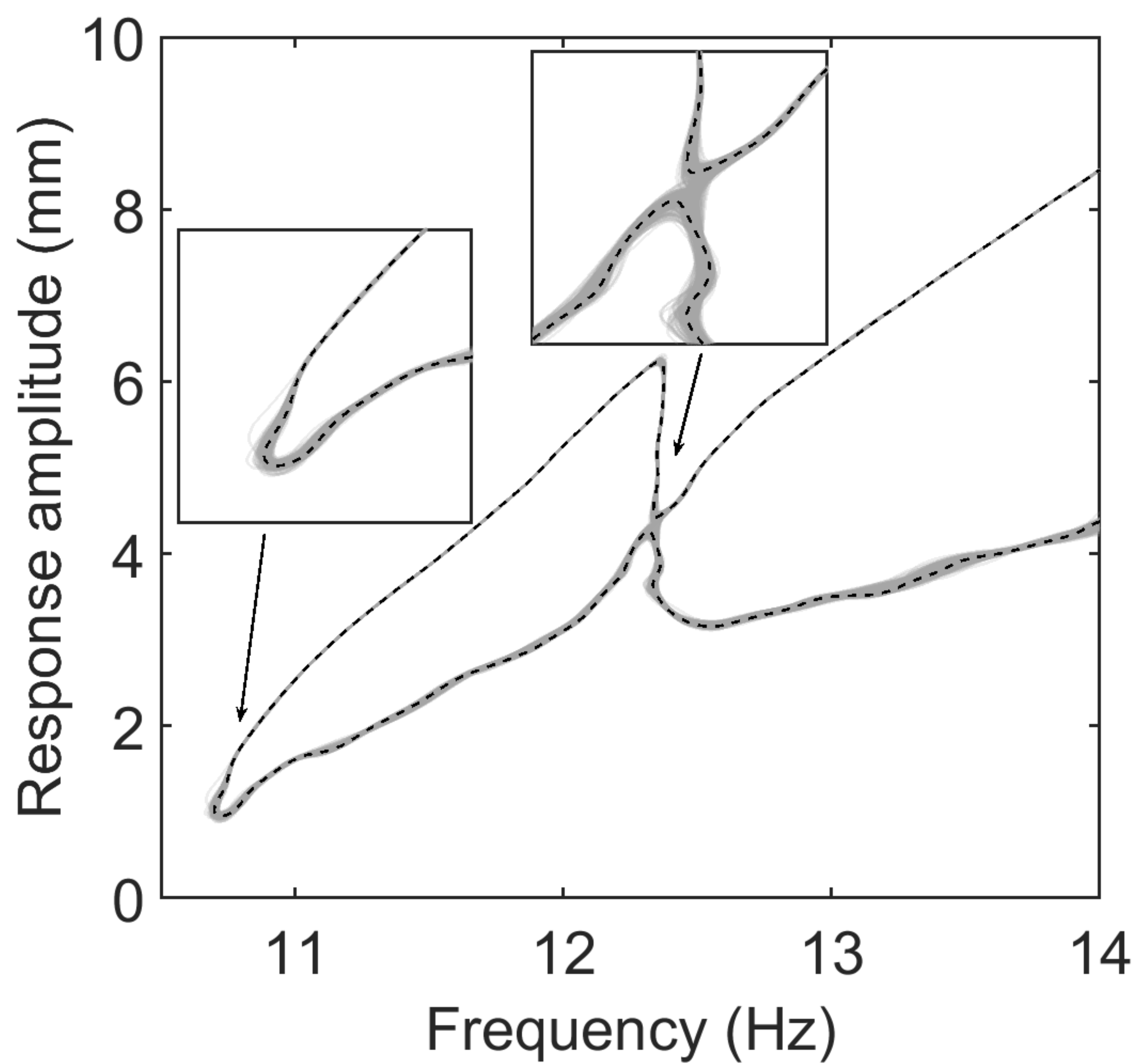}} &
\subfloat[]{\label{uq2}\includegraphics[width=0.45\textwidth]{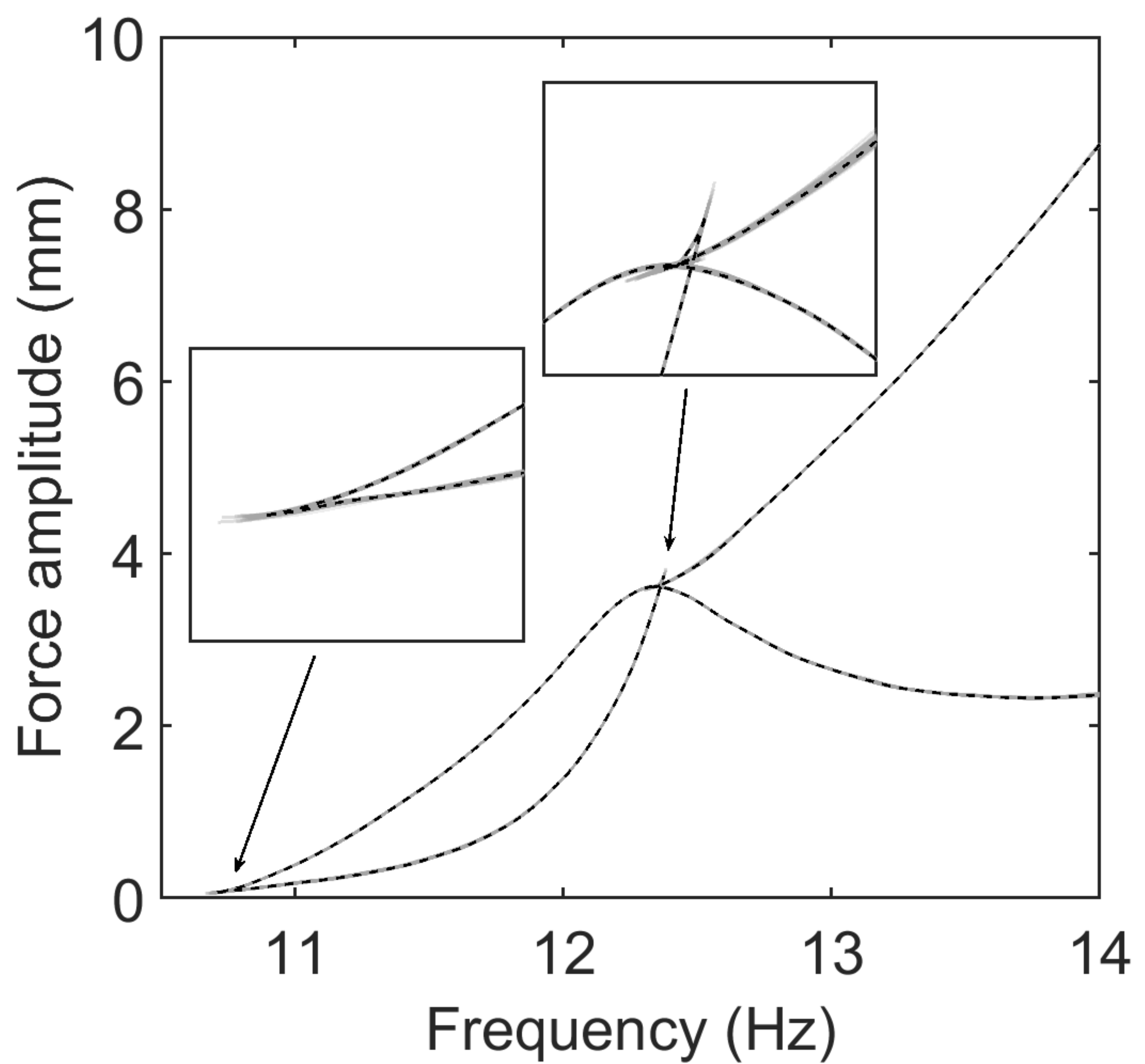}}
\end{tabular*}
\caption{(\textcolor{black}{$\boldsymbol{--}$}) fold curve obtained by post-processing all the data points collected using the online continuation algorithm. (\textcolor{lgray}{$\boldsymbol{-}$}) 300 fold curves obtained by post-processing some of the data points measured using online continuation. For each curve, 10\% of the data points are randomly selected and removed from the data set. Each curve is also computed with hyper-parameters values obtained via the maximization of the marginal likelihood for the particular data set considered.}
\label{fig:uncertainty}
\end{figure*}

The resulting fold curves are compared in Figure~\ref{fig:vs_ref}. The two curves (\textcolor{dgreen}{$\boldsymbol{-}$}, \textcolor{dblue}{$\boldsymbol{-}$}) obtained by post-processing the response surface are almost identical, which shows that the difference between the two sets of hyper-parameters does not play a significant role here. These curves are found to oscillate around the fold curve measured online (\textcolor{black}{$\boldsymbol{- \bullet -}$}). These oscillations are artifact created by the regression method due to the lack of data points (in particular at resonance between successive $S$-curves). These oscillations are not present in the curve obtained directly during the experiment due to the tailored collection of suitable data points during the continuation. The blue curve is also found to reproduce the swallowtail catastrophe observed on the data collected online. The post-processing applied to the data points of the response surface could have also been applied to the data points collected by the online regression-based continuation algorithm. The fold curve obtained in this way (\textcolor{black}{$\boldsymbol{--}$} in Figure~\ref{fig:uncertainty}) is a smooth and continuous curve that is no longer affected by the `jump' observed during the experiment and that also reproduces the swallowtail catastrophe. The curve is also free from regression artifacts as collected data points were suitably chosen.\\

The collection of data points not only at the dynamical feature of interest but also around it enables us to analyse the sensitivity of this feature to experimentally collected data points and distinguish between features that are robust (such as the swallowtail) and those that are genuinely uncertain (fold connection near the singularity). To visualize this uncertainty, 300 fold curves were numerically computed for different data sets. Each set is based on all the data points collected using the online algorithm but with 10\% of the points, selected randomly, removed. The hyper-parameters of the GPR model were marginalized for each data set. Fold curves (\textcolor{dgray}{$\boldsymbol{-}$}) are shown in Figure~\ref{fig:uncertainty}. Regions of larger result variability, i.e. where the overall bifurcation curve appears thicker, are located on the lower part of the curve (Figure~\ref{fig:uncertainty}(a)). This is expected as this region of the response surface is significantly less curved than the resonance region and hence more affected by measurement uncertainties. Large variability is also noticeable at the cusp as well as at the modal interaction region. In the latter, the distinction between the upper and lower part of the fold curve can disappear as some fold curves are found to cross and have a branch point. Other continuation runs are also found to lead to disjoint fold curves --- one before and one after the modal interaction. This important variability in the results highlights the ``flat'' nature of the response surface in that region and hence the increased difficulty in inferring the presence of folds (inset in Figure~\ref{fig:uncertainty}(a)). In general, the presence of uncertainty may make it difficult to draw strong conclusions on the precise geometry of the fold curve in such a region. Figure~\ref{fig:uncertainty}(b) shows the sensitivity of the results in terms of the bifurcations parameters forcing frequency and amplitude. As observed experimentally, the swallowtail catastrophe appears to be a robust feature of the system's dynamics as it is present for almost all data sets.

\section{Conclusions}\label{sec:conclusion}
Control-based continuation is a general and systematic method to probe the dynamics of nonlinear experiments. In this paper, CBC is combined with a novel continuation algorithm that is robust to experimental noise and enables the tracking of geometric features of the response surface. The online regression-based CBC uses Gaussian process regression to create local models of the response surface. These models are smooth (noise-free) which enables the use of standard numerical continuation algorithms. An novel aspect of the proposed method is a separate data collection step. This step selects input values maximizing information about the dynamic feature of interest. This approach is similar in principle to an experimental design problem in statistics and could a priori be applied to other problems such as nonlinear parameter estimation and model selection.\\

The online regression-based algorithm was experimentally demonstrated on a nonlinear structure with harmonically-coupled modes by tracking fold points in the structure response to harmonic excitation. The algorithm was able to address the complex geometry of the fold curve arising from the presence of a modal interaction. However, the close proximity between two portions of the fold curve lead, within uncertainties in the experiment, to a branch point singularity which could not be addressed by the current continuation algorithm.\\

The online regression-based algorithm presented in this paper is very general and could be exploited to capture other types of dynamic features such as backbone curves (nonlinear normal modes) and nonlinear frequency response curves. The regression algorithm used here relies however on a unique parameterisation of the solution curve in terms of the continuation variables. If this assumption is not satisfied, further improvements of the method to consider local coordinates for the regression, or a larger set of continuation variables might be necessary.\\

\section*{Data Statement}
Experimental data collected in this study are available at [\textit{DOI to be inserted at proofing}].

\section*{Acknowledgements}
L.R. has received funding from the Royal Academy of Engineering, fellowship RF1516/15/11. J.S. has been supported in part by funding from the European Union Horizon 2020 research and innovation programme for the ITN CRITICS under Grant Agreement Number 643073, and by EPSRC Grants EP/N023544/1 and EP/N014391/1. D.A.W.B. is funded by the EPSRC grant EP/K032738/1 and S.A.N. by the EPSRC fellowship EP/K005375/1. We gratefully acknowledge the financial support of the Royal Academy of Engineering and the EPSRC.

\bibliographystyle{unsrt}
\bibliography{mybib}

\end{document}